\magnification = 1200
\parindent=0cm

\input amstex

\define\nl{\bigskip\item{}}
\define\snl{\medskip\item{}}
\loadmsbm
\input amssym.tex
\input amssym.def

\define\ot{\otimes}

\define\jo{\iota}

\define\inspr #1{\parindent=20pt\bigskip\bf\item{#1}}
\define\iinspr #1{\parindent=27pt\bigskip\bf\item{#1}}
\define\einspr{\parindent=0cm\bigskip}
%\define\dotplus{+}

%\input qd0.tex    %  (Abstract and 0. Introduction)

\centerline{\bf Quasi-Discrete Locally Compact Quantum Groups \rm ($^\ast$)}
\nl
\centerline{\it A. Van Daele \rm ($^{\ast\ast}$)}
\nl\nl
\bf Abstract  \rm
\nl
Let $A$ be a  C$^\ast$-algebra.
Let $A\ot A$ be the minimal C$^*$-tensor product of $A$
with itself and
let $M(A\ot A)$ be the multiplier algebra of $A\ot A$.
 A comultiplication on $A$ is a non-degenerate
$^\ast$-homomorphism $\Delta : A \to M(A \otimes A)$
satisfying the coassociativity law
$(\Delta \otimes \iota)\Delta = (\iota \otimes \Delta)\Delta$
where $\iota$ is the identity map and where $\Delta\ot\iota$ and
$\iota\ot\Delta$ are the unique extensions to
 $M(A\ot A)$ of the obvious maps on $A\ot A$. We think of  a pair
$(A,\Delta)$  as a `locally compact quantum semi-group'.
\snl
When these notes where written, in 1993, it was not at all clear  
what the extra conditions on $\Delta$ should be for $(A,\Delta)$  to
be a `locally compact quantum group'. This only became clear in 1999 thanks to the work of Kustermans and Vaes. 
In the compact case however, that is when  $A$ has an identity, rather natural
conditions can be formulated and so there was a good notion of a
 'compact quantum group' already at the time these notes have been written.
These compact quantum groups have been studied by
Woronowicz. 
\snl
In these notes, we consider another class of locally compact quantum
groups. We assume the existence of a non-zero element $h$ in $A$ such
that $\Delta(a)(1\ot h)=a\ot h$ for all $a\in A$. With some extra, but
also natural conditions, the element $h$ is unique.
We speak of a {\it quasi-discrete locally compact
quantum group}.
We also discuss the discrete case and we show that, in that case, there
exists such an element $h$. So, the quasi-discrete case is, at least in
principle, more general than the discrete case. Later however, it has been shown by Kustermans that a quasi-discrete locally compact quantum group has to be a discrete quantum group.
\snl
We prove the existence of the Haar measure, the regular representation,
the fundamental unitary that satifies the Pentagon equation and we
obtain the reduced dual. 
\snl
These notes have not been published. Nevertheless, some of the results and techniques seem to be useful and in recent work, we came across similar settings. Therefore, we have decided to publish these notes in the archive. We have added some comments at the end of the introduction, and also updated the reference list. But apart from these minor changes, the notes are still as they were written in 1993. 
\medskip
\hskip 8 cm August 2004
\bigskip
\noindent
\hrule
\medskip
\parindent 0.8 cm
\item{($^\ast$)} Lecture Notes K.U.Leuven (June 1993)
\smallskip
\item{($^{\ast \ast}$)} Address: Department of Mathematics, K.U.Leuven, Celestijnenlaan
200B, B-3001 Heverlee (Belgium).

\newpage
\parindent 0 cm
\bf 0. Introduction \rm
\nl
Let us first recall the main ideas behind the
C$^\ast$-approach to quantum groups.
\snl
Let $G$ be a locally compact group and
 let $A$ be the C$^\ast$-algebra $C_0(G)$ of
continuous complex functions on $G$ tending to $0$ at infinity.
We identify the C$^\ast$-tensor product $A \overline{\otimes} A$
with $C_0 (G \times G)$ and the multiplier algebra
$M(A \overline{\otimes} A)$ with the
C$^\ast$-algebra $C_b(G \times G)$ of bounded continuous complex functions
on $G \times G$.  Then, the  product in $G$ gives rise to a
$^\ast$-homomorphism $\Delta : A \rightarrow M(A \overline{\otimes} A)$,
called the comultiplication. It
is defined by $(\Delta (f))(s,t) = f(st)$ when $f \in C_0(G)$ and $s,t \in G$.
It is well known that the topological structure of $G$ is completely determined
by the C$^\ast$-algebra structure of $A$.
Moreover, it is easy to see that the map $\Delta$ completely determines
the group structure.
\snl
The idea of the C$^\ast$-algebra approach to quantum groups is now very
natural. The abelian C$^*$-algebra $A$ is replaced by any C$^*$-algebra.
So a `locally compact quantum group' is a pair $(A,\Delta)$ where now $A$
is any C$^\ast$-algebra and $\Delta$ a $^\ast$-homomorphism of $A$ into $M(A
\overline{\otimes} A)$ satisfying certain properties.
\snl
It is clear that we need extra conditions on $\Delta$. Any abelian
C$^*$-algebra $A$ has the form $C_0(G)$ for some locally compact space
$G$ but not every $^*$-homomorphism $\Delta : A\to M(A\overline\ot A)$
will come from a group structure on $G$ as above.
\snl
Some conditions on $\Delta $ are quite natural. First, $\Delta$ has to
be a non-degenerate $^*$-homomor\-phism in the sense of [29]. This means
that $\Delta (A)(A\ot A)$ is dense in $A\overline\ot A$. In the abelian
case, i.e.\ when $A\cong C_0(G)$, with $G$ a locally compact space, this
will guarantee that $\Delta$ comes from a product  on $G$. The
coassociativity law $(\Delta \ot\iota)\Delta=(\iota\ot\Delta)\Delta$ will
give that this procuct is associative. Therefore, a C$^*$-algebra $A$
with a non-degenerate $^*$-homomorphism $\Delta : A\to M(A\overline\ot
A)$ satsifying coassociativity can be thought of as a locally compact
quantum semi-group.
\snl
What extra conditions are needed in order to have a locally compact
quantum group? The answer to this question is not yet clear. We want
the conditions to be such that most of the theory of locally compact
groups can be generalized to the quantum groups. We like to prove the
existence of a unique Haar measure, to have
a nice representation theory, ...
We like to define a dual pair $(\hat A,\hat\Delta)$ that is again a
locally compact quantum group (corresponding to the dual group $\hat G$
of $G$ if $A=C_0(G)$). We would like to have that $(A,\Delta)$ is again
the dual pair associated to $(\hat A,\hat\Delta)$, ...
\snl
In the compact case, that is when $A$ has an identity, most of this
program has been carried out by Woronowicz. He imposed rather natural
conditions on $\Delta$ to have a 'compact quantum group'. He proved the
existence of a Haar measure in this case. He also developed the
representation theory, he obtained the dual, ... [30, 31].
\snl
On the other hand, the discrete case has been studied by Podle\'s and
Woronowicz in [17]. In their work, the discrete quantum groups are the
duals of the compact ones. The discrete case is also studied by Effros
and Ruan [6]. But there, the approach is more algebraic. It does not
really fit into the C$^*$-algebra framework.
\snl
In this paper, we study the discrete case independently. In the
philosophy of the C$^*$-algebra approach to quantum groups, discreteness
is a property of the topology of the group and therefore, in the quantum
group case, we must translate it into a property of the underlying
C$^*$-algebra. This is not so difficult, we discuss this in section 2 of
this paper.
\snl
However, there is also the following simple observation.
If $G$ is a discrete
group and $h$ is the function $\delta_e$ defined as $1$ on the identity
$e$ and $0$ elsewhere, then $h$ is an element in $C_0(G)$ such that
$(\Delta (f))(1\ot h)= f\ot h$ for all $f\in C_0(G)$.
Indeed, when $s,t\in G$ we have
$$(\Delta (f))(1\ot h)(s,t)=f(st)\delta_e(t)=f(s)\delta_e(t)=(f\ot
h)(s,t).$$
\snl
We can also take the
existence of such an element as an axiom. So we introduce the
notion of a quasi-discrete quantum group as a pair $(A,\Delta)$ of a
C$^*$-algebra $A$ and a comultiplication $\Delta$ such that there is a
non-zero element $h\in A$  with the property that $\Delta(a)(1\ot
h)=a\ot h$ for all $a\in A$. Of course also here we need some extra
conditions to distinguish from the semi-group case.
\snl
It turns out that in the discrete case, such an element $h$
automatically exists. So, a discrete quantum group is also
quasi-discrete. It is not clear however if the second class is really
bigger that the first one. We have no examples to show this. Certainly,
the techniques that we use are different from the ones that are normally
used in the discrete case (where $A$ is a direct sum of full matrix
algebras).  So, even though the two classes might turn out to be the
same, still we have achieved two goals. On the one hand, we have a
self-contained treatment of the discrete quantum groups. On the other
hand, we have used techniques that will probably be useful to
develop  the general locally compact quantum groups.
\nl
The paper is organized as follows. In {\it Section} 1, we try to find the
natural extra conditions on $\Delta$. We focus on the abelian case here.
In {\it Section} 2 we look at the (proper) discrete case and we
compare this with the work of Effros and Ruan. In {\it Section} 3 we
give the precise definition of the quasi-discrete locally compact
quantum groups.
\snl
In the following sections we develop the theory for our quasi-disrete
locally compact quantum groups. The main point is the construction of
the antipode $S$ as a linear operator on $A$.  This is done in {\it Section} 4.
We must mention here that this is another main difference with the
approach of Effros and Ruan. The existence of an antipode is one of
their axioms. In the C$^*$-algebra approach to quantum groups however,
it is not so natural to assume the existence of the antipode. It is
often a basic difficulty that the antipode is an unbounded,
anti-homomorphism. This is also one of the main differences with the
earlier von Neumann algebra approach (the Kac algebras) where the
antipode was assumed to be bounded, in fact with square one. This older
theory was no longer satisfactory since the more recent interesting
examples were discovered.
\snl
We obtain some nice properties of the antipode that turn
out to be very useful for the next sections. In {\it Section} 5 we use it to
construct the Haar measure and in {\it Section} 6 to obtain the regular
representation. In {\it Section} 7, we construct the fundamental unitary that
satisfies the Pentagon equation. There, we also obtain the reduced dual.
\snl
In the discrete case, the Haar weight is semi-finite, but in the
more general quasi-discrete case, it is not. Perhaps this is not so
surprising. After all, in the compact case, the Haar measure need not be
faithful. And in a way, these properties are dual to each other. The
existence of compact quantum groups with a non-faithful Haar measure therefore
may indicate that there are quasi-discrete quantum groups that are not
discrete.
\nl
In these notes, we will mostly work with separable C$^*$-algebras for technical
convenience. We believe though, that this is not essential and that all
the arguments can be formulated also in the non-separable case.
\nl
\it June 1993 \rm
\nl\nl
\bf Note \rm
\nl
These lecture notes have been written in 1993. At that time, there was a clear notion of a 'compact quantum group' (see [30, 31]) and about simultaneously with the appearance of these notes, a notion of 'discrete quantum groups' was developed by Effros and Ruan (see [6]). The notion of a locally compact quantum group was fully established in 1999 by Kustermans and Vaes (see [11], [12] and [13]).
\snl
In these notes, a certain class of 'locally compact quantum groups' is developed, the so-called quasi-discrete locally compact quantum groups. Shortly after this work was done, it was shown by Kustermans [10], that a quasi-discrete locally compact quantum group was actually a discrete quantum group (as introduced by Effros and Ruan in [6] and later by myself in [26]). For this reason, we decided not to publish these notes. Also the work by Kustermans has not been published. 
\snl
Still, we believe that this work is of some interest. Recently, we have been involved in some research where a similar setting arose (although purely algebraic), see [14] and [28]. This is the reason why we have decided to publish this paper on the net. We have not made changes to the original version of 1993. Except for the abstract, this note and the references (which have been updated), the paper is the same as in 1993. This has to be taken into account by the reader.  
\nl
\it August 2004 \rm
\nl\nl

\bf Acknowledgements \rm
\nl
Most of this work was done while we were visiting the University of Oslo in the
summer of 1992.  We were greatly inspired by  many discussions we had there
with E. Alfsen on the Haar measure (for locally compact groups) and we
are  grateful to him, both for these discussions and for his hospitality
while we were staying in Oslo.
\nl\nl

\bf 1. Locally compact quantum groups \rm
\nl
Consider an abelian C$^\ast$-algebra $A$ and a $^\ast$-homomorphism
$\Delta : A \rightarrow M(A \overline{\otimes} A)$. We know that
$A$ has the form $C_0(X)$ where $X$ is a locally compact space.
We will look for
conditions on $\Delta$ to have that it is of the form
$(\Delta f)(s,t)=f(st)$ for a multiplication on $X$ that makes $X$ into a
locally compact group.
\snl
We will not prove very deep results in this section.  It will
mainly serve as a motivation and we will use the abelian case to illustrate
some aspects in the general case.  We are
interested in finding a set of axioms that will also work in the
non-abelian case.
So, when dealing here with the abelian case, we must try to avoid the
abelianness as much as possible and certainly in the formulation of the
conditions.
\nl
So, assume that  $A= C_0(X)$ where $X$ is a locally compact space. Then,
the $^\ast$-homomor\-phism $\Delta : A
\rightarrow  M(A \overline{\otimes} A)$ gives a mapping from $C_0(X)$ into $C_b(X
\times X)$. The following regularity condition is sufficient to have that
$\Delta$ is given by a product law on $X$.

\inspr{1.1} Proposition. \rm Assume that $\Delta (A) (A \otimes A)$ is
dense in $A \overline{\otimes} A$.  Then there is a continuous map $(s,t)
\in X \times X \rightarrow st \in X$ such that $(\Delta f)(s,t) = f(st)$ when $f
\in A$ and $s,t \in X$.
\snl\bf Proof : \rm  Choose two elements $s,t \in X$ and consider the map $f \rightarrow
(\Delta f)(s,t)$.  It is clear that this is a $^\ast$-homomorphism from $A$ to
$\Bbb C$.  If $(\Delta f)(s,t) = 0$ for all $f$, then $((\Delta f)(g \otimes
h))(s,t) = 0$ for all $f,g,h \in A$.  By the density condition, $f(s)g(t) = (f
\otimes g)(s,t) = 0$ for all $f,g \in A$.  This is impossible.  Hence the above
$^\ast$-homomorphism is non-zero.  Therefore, there is an element in $X$, denoted
by $st$, such that $(\Delta f)(s,t) = f(st)$ for all $f \in A$.
\snl
This proves the existence of the product.  The continuity follows from the fact
that the map $(s,t) \rightarrow f(st)$ is continuous for all $f \in A$.
\einspr

A $^\ast$-homomorphism $\Delta : A \rightarrow M(A \overline{\otimes} A)$
satisfying the condition that $\Delta(A)(A \otimes A)$ is
dense in $A
\overline{\otimes} A$ is called non-degenerate.  It is a morphism from $A$
to $A \overline{\otimes} A$ in the sense of [29].
\nl
We now want the product in $X$ to be associative.  For this we need the
coassociativity of $\Delta$, as formulated in the following proposition.

\inspr{1.2} Proposition \rm
Assume that $\Delta$ is non-degenerate.  Consider the
$^\ast$-homomorphisms $\Delta \otimes \iota$ and $\iota \otimes \Delta$ on $A
\otimes A$ where $\iota$ is the identity map.  They have unique
extensions to $^\ast$-homomorphisms on $M(A \overline{\otimes} A)$.  We
denote these extensions still by $\Delta \otimes \iota$ and $\iota \otimes
\Delta$.  Assume that $(\Delta \otimes \iota)\Delta = (\iota \otimes
\Delta)\Delta$.  Then, the product on $X$, defined by $\Delta$ as in proposition
1.1, is associative.

\snl\bf Proof : \rm
It is not so hard to show that the $^\ast$-homomorphisms $\Delta \otimes
\iota$ and $\iota \otimes \Delta$ have unique extensions.  This can be done,
also when $A$ is non-abelian.  In the abelian case, it is in fact obvious.  It is
clear that the extensions are given by
$((\Delta \otimes \iota)f)(r,s,t) = f(rs,t)$
and $((\iota \otimes \Delta)f)(r,s,t) = f(r,st)$
whenever $f \in C_b (X \times X)$ and $r,s,t \in X$.  Then it is also clear that
the condition $(\Delta \otimes \iota)\Delta = (\iota \otimes \Delta)\Delta$
gives the associativity of the product.  Indeed
$$\align ((\Delta \otimes \iota)\Delta f)(r,s,t) &
= (\Delta f)(rs,t) = f((rs)t)\\
((\iota \otimes \Delta)\Delta f)(r,s,t) &= (\Delta f)(r,st) =
f(r(st)).\endalign$$
\einspr

These two conditions bring us to the following definitions. Here, $A$ is
again any C$^*$-algebra.

\inspr{1.3} Definition \rm
Let $A$ be a C$^\ast$-algebra and $\Delta : A \rightarrow
M(A \overline{\otimes} A)$ a $^\ast$-homomorphism, where $A
\overline{\otimes} A$ is some C$^\ast$-tensor product of $A$ with itself.  Assume that
$\Delta$ is non-degenerate, i.e. that $\Delta(A)(A \otimes A)$ is
dense in $A \overline{\otimes} A$.  Also assume that $\Delta$ is
coassociative, i.e. that $(\Delta \otimes \iota)\Delta = (\iota \otimes
\Delta)\Delta$.  Then $\Delta$ is called a comultiplication on $A$.  We will also
call the pair $(A,\Delta)$ a locally compact quantum semi-group.
\einspr

We now proceed by looking for conditions that are sufficient, in the abelian
case, to have a group.  So, in what follows, $A$ is again an abelian
C$^\ast$-algebra $C_0(X)$ and $\Delta$ is a comultiplication on $A$.
\nl
The following conditions are quite natural.

\inspr{1.4} Proposition \rm
Assume that $\Delta(A)(1 \otimes A)$ and $\Delta(A)(A \otimes
1)$ are dense subspaces of $A \overline{\otimes} A$.  Then the product in $X$,
associated with $\Delta$, has the cancellation law, i.e. if $st = rt$ then $s=r$
and if $ts = tr$ then $s=r$.

\snl\bf Proof : \rm
Suppose that $r,s,t \in X$ and that $st = rt$.  Then, for all $f,g \in
A$ we have
$$\align(\Delta(f)(1 \otimes g))(s,t) &= f(st)g(t)\\
(\Delta(f)(1 \otimes g))(r,t) &= f(rt)g(t)\endalign$$
so that
$$(\Delta (f)(1 \otimes g))(s,t) = (\Delta (f)(1 \otimes g))(r,t).$$
By the density of $\Delta(A)(1 \otimes A)$ in $A \overline{\otimes} A$, we
have $f(s)g(t) = f(r)g(t)$ for all $f,g \in A$.  Hence $r=s$.
\snl
Similarly, if $ts = tr$, the density of $\Delta(A)(A \otimes 1)$ will give $s =
r$.
\einspr

We want to make a remark here about the relation of these conditions and the
previous ones.  First it is clear that the density of $\Delta(A)(1 \otimes A)$
in $A \overline{\otimes} A$  automatically gives the density of $\Delta(A)(A
\otimes A)$ in $A \overline{\otimes} A$, and similarly for $\Delta(A)(A
\otimes 1)$.  A second remark is that one can now formulate
 the coassociativity
rule, as
$$(a \otimes 1)(\Delta \otimes \iota)(\Delta(b)(1 \otimes c)) = (\iota \otimes
\Delta)((a \otimes 1)\Delta(b)) (1 \otimes c),$$
that is, without having the need to consider the extensions of $\Delta \otimes
\iota$ and $\iota \otimes \Delta$ to the multiplier algebra (see e.g.
[24]).
\nl
Now, it is known that a compact semi-group with cancellation is actually
a group, see [8]. We give here, for completeness, a simple proof of this
result.

\inspr{1.5} Propostion \rm
Let $G$ be a compact semi-group with cancellation, then $G$ is a group.

\snl\bf Proof : \rm
Take any element $s\in G$. Consider the closed semi-group $H$ generated
by $s$. Consider the family of all closed ideals in $H$. The
intersection of two ideals is non-empty and again an ideal.
Because $H$ is compact,
the intersection of all ideals is an ideal $I$. It is the minimal ideal in
$H$. Because $H$ is abelian, $tI$ is again an ideal for all $t$ in $H$.
It is contained in $I$ and by minimality, it is equal to $I$. So, $tI=I$
for all $t\in H$.
\snl
Consider any element $t\in I$ and choose $e\in I$ scuh that
$te=t$.
Multiply with any element $r\in G$ to the right and cancel $t$ to get $er=r$.
Then multiply with any element $t\in G$ to the left and cancel $r$ to get
$te=t$. So, $e$ is the identity in $G$.
\snl
Now, from $se=s$ and the fact that $e\in I$ and that $I$ is an ideal, we
must have $s\in I$. Because  $sI=I$ and $e\in I$, $s$ has an inverse. So $G$ is
a group.
\einspr

So, if $A$ is any abelian  C$^\ast$-algebra with an
identity and if $\Delta$ is a comultiplication on $A$ such that $\Delta(A)(1
\otimes A)$ and $\Delta(A)(A \otimes 1)$ are dense subspaces of $A
\overline{\otimes}A$, then $A \cong C_0(G)$ where $G$ is a compact group and
$\Delta$ is given by the multiplication on $G$ as before.
\nl
So we are lead naturally  to the following definition, due to Woronowicz
[31] :

\inspr{1.6} Definition \rm
If $A$ is any C$^\ast$-algebra with identity and $\Delta$ a
comultiplication on $A$, such that the sets $\Delta(A)(1 \otimes A)$ and
$\Delta(A)(A \otimes 1)$ are dense in $A \overline{\otimes} A$, then
$(A,\Delta)$ is called a compact quantum group.
\einspr

Woronowicz proved the existence of the Haar measure for such compact quantum
groups when $A$ admits a faithful state, see [31]. We adapted his proof to
obtain the Haar measure without this restriction on $A$, see [25].
For compact quantum groups, the Haar measure is unique but,
in general, it need not be faithful.
Nevertheless, much of the theory of compact groups extends to the compact
quantum groups as defined above, see [30, 31].
\nl
In the non-compact case, the above conditions are not sufficient to have a
group.  We illustrate this with some simple examples.

\inspr{1.7} Examples \rm
i) Consider the (discrete) semi-group $G = \Bbb N \cup \{ \infty \}$ (with
$n + \infty = \infty + n = \infty$ for all $n \in \Bbb N$).  Let $\delta_p$
denote the function that is $1$ on $p$ and $0$ elsewhere. We have that
$(\Delta \delta_\infty)(n,m) = \delta_\infty (n+m).$
Now $n+m = \infty$ if and only if $n = \infty$ or $m = \infty$.  So we get
$$\Delta (\delta_\infty) = \delta_\infty \otimes 1 + 1 \otimes \delta_\infty -
\delta_\infty \otimes \delta_\infty.$$
In particular $\Delta(\delta_\infty)(\delta_\infty \otimes 1) = \delta_\infty
\otimes 1$ and we see that in this example $\Delta(A)(A \otimes 1)$ is not a
subspace of $A \overline{\otimes} A$.  Of course, also here, we do not have
cancellation.
\nl
ii) Let $G = \Bbb N$.  For any pair $(n,m)$ in $\Bbb N$ we have
$$\Delta (\delta_n) (1 \otimes \delta_m) = \cases \delta_{n-m} \otimes \delta_m
& \text{if $n \geq m$} \\
0 & \text{if $n < m$}. \endcases
$$
In particular, $\Delta(A)(1 \otimes A)$ is a dense subspace of $A
\overline{\otimes} A$ here.  Similarly of course $\Delta(A) (A \otimes 1)$ is a dense
subspace of $A
\overline{\otimes} A$.  Indeed, $\Bbb N$ has the cancellation
property.  Still, $\Bbb N$ is not a group.
\einspr

In the last example, we see that there are elements $a,b \in A$ such that
$\Delta(a)(1 \otimes b) = 0$ with $a \ne 0$ and $b \ne 0$.  This cannot happen
in the group case : if $f(st)g(t) = 0$ for all $s,t \in G$, then $f(r)g(t) = 0$
for all $r,t \in G$ and hence $f = 0$ or $g = 0$.
\nl
It turns out that this is precisely the extra condition that we need in general.

\inspr{1.8} Theorem \rm
Let $A$ be an abelian C$^\ast$-algebra with a comultiplication
$\Delta$ such that $\Delta(A)(1 \otimes A)$ and $\Delta(A)(A \otimes 1)$ are
dense subspaces of $A \overline{\otimes} A$.  Also assume that $\Delta(a)(1
\otimes b)  =0$ with $a,b \in A$ implies $a = 0$ or $b = 0$.  Then $A \cong
C_0(G)$ where $G$ is a locally compact group and $\Delta$ is given by the
multiplication on $G$ as before.

\snl\bf Proof : \rm
We know already that $A \cong C_0(G)$ and that $G$ is a semi-group with the
cancellation property.  Let us now see what the last injectivity condition
means in terms of $G$.
\snl
Choose any non-empty open subset $V$ of $G$.  We claim that the set $\{ r s
\mid r \in
G, s \in V \}$ is dense in $G$.  Indeed, suppose that it is not dense.  Then
there is an element $f \in C_0(G)$ such that $f \ne 0$ but $f(rs) = 0$ for all
$r \in G$ and all $s \in V$.  Choose an element $g \in C_0(G)$ such that $g \ne
0$ and such that $g$ has support in $V$.  Then $f(rs)g(s) = 0$ for all $r,s \in
G$.  Hence $(\Delta f) (1 \otimes g) = 0$.  This would imply $f = 0$ or $g = 0$ by
assumption.  So we get a contradiction.
\snl
Now, let $s,t \in G$.  For any pair $V,W$ of open neighbourhoods of $s$ and $t$
respectively, we have, by the above property, a point $r$ in $G$ such that $rs =
t$.  Consider the pairs $(V,W)$ as an index set $I$, ordered by the opposite
inclusion.  Then we find nets $(r_\alpha), (s_\alpha), (t_\alpha)$ such that
$r_\alpha s_\alpha = t_\alpha$ for all $\alpha \in I$ and such that $s_\alpha
\rightarrow s$ and $t_\alpha \rightarrow t$.
\snl
For every $f,g \in A$ we have
$$((\Delta f)(1 \otimes g))(r_\alpha,s_\alpha) = f(r_\alpha s_\alpha)g(s_\alpha)
= f(t_\alpha)g(s_\alpha)$$
and this converges to $f(t)g(s)$.  Choose $f$ and $g$ such that $f(t)g(s) \ne
0$.  Let $K_1,K_2$ be compact sets in $G$ such that $(\Delta f)(1 \otimes g)$ is
strictly smaller than $\frac 12 | f(t) g(s) |$ outside $K_1 \times K_2$.
Then, for $\alpha$ large enough, we will have that
$|((\Delta f)(1\otimes g))(r_\alpha,s_\alpha)|\geq \frac12|f(t)g(s)|$ and
so $r_\alpha\in K_1$ and $s_\alpha\in K_2$. Then,
we can find a subnet of $(r_\alpha)$ that converges to a point $r$.  So in the
limit, we will get
$rs = t$.
\snl
So for every two elements $s,t \in G$ we find $r \in G$ such that $rs = t$.
Given $s \in G$ we get $e \in G$ such that $es = s$.  As in the proof of 1.5, we
see, using cancellation that $e$ is an identity.  Again, given $s$ we find $r$
such that $rs = e$. Multiplying with $s$ on the left and cancelling $s$ on the
right, we get $sr = e$.  So every element has an inverse.
\snl
We also must show that the maps $s\to s^{-1}$ is continuous. For all
$f,g\in A$ we have that
$((\Delta f)(1\ot g))(s^{-1},s)=f(s^{-1}s)g(s)=f(e)g(s)$ so that
$s\to ((\Delta f)(1\ot g))(s^{-1},s)$
is continuous. By density of $\Delta(A)(1\ot A)$ in $A\overline\ot A$,
still $s\to f(s^{-1})g(s)$ will be continuous for all $f,g\in A$.
Hence,  $s\to s^{-1}$ is continuous.
\einspr

It is clear that the condition $\Delta (a)(1\ot b)=0$ implies $a=0$ or
$b=0$  can be replaced by the similar condition $\Delta (a)(b\ot 1)=0$
implies $a=0$ or $b=0$. But there is more. If we examine the proof of the
previous theorem, we see that we also have the following result.

\inspr{1.9} Theorem \rm
Let $A$ be an abelian C$^\ast$-algebra with a comultiplication
$\Delta$ such that $\Delta(A)(1 \otimes A)$ and $\Delta(A)(A\ot 1)$
are subspaces of $A \overline{\otimes} A$.  Assume further that
$$\align\Delta(a)(1 \otimes b) &= 0 \Rightarrow a \otimes b = 0\\
\Delta(a)(b \otimes 1) &= 0 \Rightarrow a \otimes b = 0.\endalign$$
Then $A \cong C_0(G)$ where $G$ is a locally compact group and $\Delta$ is given
by the multiplication on $G$.

\snl\bf Proof : \rm
The injectivity rules give that, for any pair $s,t \in G$,
there exists elements $u$ and $v$ such that $su = t$ and $vs = t$.
This was shown in the previous theorem. Remark that we need that
  $\Delta(A)(1 \otimes A)$ and $\Delta(A)(A\ot 1)$
are subsets of $A \overline{\otimes} A$. We do not need the density of
these spaces. So $tG=Gt=G$ for all $t \in G$. Then $G$ is a group.
\snl
The continuity of $s \to s^{-1}$ is proven as in the previous case.
\einspr

We finish this section with some important remarks.

\iinspr{1.10} Remarks \rm
i) Consider the linear map $T : A \otimes A \rightarrow A
\overline{\otimes} A$ defined by $T(a \otimes b) = \Delta(a)(1 \otimes b)$.
In the group case $A = C_0(G)$, this map extends to an isomorphism on $A
\overline{\otimes} A$ given by $(Tf)(s,t) = f(st,t)$.  The inverse is of
course given by $(T^{-1}f)(s,t) = f(st^{-1},t)$.
\snl
ii) We see in the above theorems that the injectivity rule $T(a \otimes b) = 0
\Rightarrow a \otimes b = 0$ implies that also the extension of $T$ remains
injective.  A similar result is also true for the compact quantum groups (see
[31]).  And we see further in this paper (see proposition 3.9 below),
that it is also true for the quasi-discrete
quantum groups.
\snl
iii) In our paper on multiplier Hopf algebras, we needed the fact that $T$ was a
bijection between the algebraic tensor product $A \otimes A$ and itself.  This,
together with the property $\Delta(A)(A \otimes 1) = A \otimes A$ was enough, in
that algebraic context, to develop the theory ([24]).
\einspr

In view of all these remarks,
for a pair $(A,\Delta)$ with any C$^*$-algebra $A$ and a
comultiplication $\Delta$ on $A$ to be a locally compact quantum group,
we will need conditions like the density of
$\Delta (A)(A\ot 1)$ and $\Delta(A)(1\ot A)$ in $A\overline\ot A$ and the
injectivity of the maps $a\ot b\to\Delta(a)(1\ot b)$ and $a\ot b\to
\Delta(a)(b\ot 1)$.
\nl\nl

\bf 2. Discrete locally compact quantum groups \rm
\nl
A discrete locally compact quantum group is, in the above philosophy, a locally
compact quantum group where the underlying space is discrete.  Therefore, we
must look at properties of $C_0(X)$ when $X$ is discrete, and see how they can
be generalized to a not necessarily abelian C$^\ast$-algebra.
\nl
If $A = C_0(X)$, then the multiplier algebra $M(A)$ of $A$ can be identified
with the C$^\ast$-algebra $C_b(X)$ of all bounded continuous complex functions
on $X$.  If $X$ is discrete, then $C_b(X)$ consists of all bounded complex
functions on $X$ and therefore coincides with the bidual of $C_0(X)$.  We can
take this property as the basis for the following definition.

\inspr{2.1} Definition \rm
A C$^\ast$-algebra $A$ such that the multiplier algebra $M(A)$
coincides with the full bidual $A^{\ast \ast}$ is called a discrete quantum
space.
\einspr

Recall that $A^{\ast \ast}$ coincides with the von Neumann algebra $A''$
when $A$ is considered in its
universal representation [16] and that then, $M(A)$
 is  the subalgebra of elements
$x$ in $A^{\prime \prime}$ such that $xa$ and $ax$ are elements in $A$ for all
$a \in A$.  So $A$ satisfies the above property if $xa \in A$ and $ax \in A$
whenever $a \in A$ and $x \in A^{\prime\prime}$.
\snl
It is not hard to see that such a C$^\ast$-algebra must be a direct sum of
components that are either full matrix algebras or
 C$^\ast$-algebras of compact
operators on infinite-dimensional Hilbert spaces  (see e.g. [15]).
\nl
So let $A = \sum_\alpha \oplus A_\alpha$ where each $A_\alpha$ is of the form
$\Cal K(\Cal H_\alpha)$, the C$^\ast$-algebra of compact operators on some
Hilbert space $\Cal H_\alpha$.
Here $\Cal H_\alpha$ can be finite-dimensional so that $\Cal K(\Cal H_\alpha)=
M_n (\Bbb C)$ where $n$ is the dimension of $\Cal H_\alpha$.  Recall that the
direct sum $\sum_{\alpha} \oplus A_\alpha$ consists of elements $(x_\alpha)_\alpha$ where
$x_\alpha \in A_\alpha$ for each $\alpha$ and $\| x_\alpha \| \rightarrow 0$ as
$\alpha \rightarrow \infty$ (the index set being considered as a discrete space
here).
\snl
We have $A \overline{\otimes} A =
\sum_{\alpha,\beta} \oplus (A_\alpha \otimes A_\beta)$ and
$M(A \overline\otimes A) = \prod_{\alpha,\beta} B(\Cal H_\alpha \otimes
\Cal H_\beta)$, where $B(\Cal H_{\alpha}\ot \Cal H_{\beta})$ denotes the
C$^*$-algebra of all bounded linear operators on
$\Cal H_{\alpha}\ot \Cal H_{\beta}$.
Recall also that the product $\prod_\alpha A_\alpha$ consists of elements
$(x_\alpha)_\alpha$ such that $x_\alpha \in A_\alpha$ for all $\alpha$ and such
that $\| x_\alpha \|$ is bounded in $\alpha$.
\snl
Now, let $\Delta : A \rightarrow M(A \overline{\otimes} A)$ be a
comultiplication on $A$.
In this case, it gives, in a
natural way, a family of $^\ast$-homomorphisms $\Delta^{\beta \gamma}_\alpha :
A_\alpha \rightarrow B(\Cal H_\beta \otimes \Cal H_\gamma)$. Then, we
have the following result (see also [6] and [17]).

\inspr{2.2}  Proposition \rm
Assume $\Delta(A)(A \otimes 1) \subseteq A \overline
\otimes  A$ and $\Delta(A)(1 \otimes A) \subseteq A \overline{\otimes} A$.
Then given $\alpha, \beta$, there exists only finitely many $\gamma$ such that
$\Delta^{\beta \gamma}_\alpha \ne 0$. Similarly, given $\alpha, \gamma$,
there exists only finitely many $\beta$ such that
$\Delta^{\beta \gamma}_\alpha \ne 0$.

\snl\bf Proof : \rm
Denote by $e_\alpha$ the identity in $A_\alpha$ for any index $\alpha$.
Fix $\alpha$ and $\beta$.  Then $\Delta(e_\alpha) (e_\beta \otimes 1)$ is a
projection in $A \overline{\otimes} A$ by assumption.  It must lie in
$\sum_\gamma \oplus (A_\beta \otimes A_\gamma)$.  Such a projection must be of
the form $(x_\gamma)_\gamma$ with $x_\gamma \in A_\beta \otimes A_\gamma$ and
$x_\gamma$ a projection for all $\gamma$.  Because $\| x_\gamma \| \rightarrow
0$ as $\gamma \rightarrow \infty$, we must have $x_\gamma = 0$ except for a
finite number of indices $\gamma$.  This proves the first statement.
Similarly for the other one.
\einspr

Now assume that all $A_{\alpha}$ are finite-dimensional. Consider the
ideal $\Cal A$ of $A$ consisting of all elements $(x_{\alpha})$ with
only finitely many $x_{\alpha}$ non-zero.  It is clear from the above
result that $\Delta (\Cal A)(\Cal A\ot 1)$ and $\Delta (\Cal A)(1\ot \Cal A)$ now are
subspaces of the algebraic tensor product  $\Cal A\ot \Cal A$
of $\Cal A$ with itself.
If the linear maps $T_1$ and $T_2$, defined on $\Cal A\ot \Cal A$ by
$T_1(a\ot
b)=\Delta(a)(1\ot b)$ and $T_2(a\ot b)=\Delta (a)(b\ot 1)$
 are bijective,
then the pair $(\Cal A,\Delta )$ is a multiplier Hopf
$^*$-algebra (in the sense of [24]).
\snl
We propose the following definition for a discrete quantum group.

\inspr{2.3} Definition \rm
A discrete quantum group is a pair $(A,\Delta)$ of a C$^*$-algebra $A$
and a comultiplication $\Delta$ on $A$ where $A$ has the form
$\sum_{\alpha} \oplus A_{\alpha}$ with each $A_{\alpha}$ a full matrix
algebra and where the linear maps $T_1$ and $T_2$, defined by
$T_1(a\ot b)=\Delta(a)(1\ot b)$ and $T_2(a\ot b)=\Delta (a)(b\ot 1)$,
give bijections of the algebraic tensor product $\Cal A\ot \Cal A$ of
the $^*$-algebra $\Cal A$ of elements $(x_{\alpha})$ in $A$ with only
finitely many $x_{\alpha}$ non-zero.
\einspr

Essentially, we get a discrete quantum group when we have a multiplier
Hopf $^*$-algebra $(\Cal A,\Delta)$ with a $^*$-algebra $\Cal A$ as
above. We also know from the theory of multiplier Hopf algebras that
there exists a unique counit $\epsilon$ and an antipode $S$ which is
invertible [24]. So we get the same objects as Effros and Ruan in [6].
The fact of working with C$^*$-algebras is no longer essential.
\snl
Given the counit $\epsilon$, we have the following result.

\inspr{2.4} Proposition \rm
Let $(A, \Delta)$ be a discrete quantum group. Let $\epsilon$ be the
counit.
Then there is a
self-adjoint projection $h$ in $A$ such that $ha = ah = \epsilon (a) h$  for all
$a \in A$.  We also have $\Delta(a)(1 \otimes h) = a \otimes h$ and $\Delta(a)(h
\otimes 1) = h \otimes a$ for all $a \in A$.

\snl \bf Proof : \rm
Since $\epsilon$ is a $^\ast$-homomorphism from $\Cal A$ to $\Bbb C$, the
kernel of $\epsilon$ is a two sided ideal in $\Cal A$ with codimension $1$.  Hence one of
the direct summands of $A$ must be one-dimensional.  Let $h$ be the identity of
this component.  Then $h$ is a self-adjoint projection in $A$ such that $ha = ah
= \epsilon (a)h$.
\snl
Clearly we have also
$$\align\Delta(a)(1 \otimes h) &= (\iota \otimes \epsilon)\Delta(a) \otimes h =
a \otimes h\\
\Delta(a)(h \otimes 1) &= h \otimes (\epsilon\otimes\iota)\Delta(a) = h
\otimes a.\endalign$$
\einspr

In the case of a discrete group, the element $h$ is of course the
function $\delta_e$ where $e$ is the identity of the group.
\snl
As we explained already in the introduction, we will take the existence
of such an element $h$ as an axiom. This will give us the notion of a
quasi-discrete quantum group.
\snl
By proposition 2.4, we will have that a discrete quantum group is
automatically quasi-discrete. So, all the results that we obtain in this
paper for quasi-discrete quantum groups are also valid for the discrete
ones. In the next sections, where appropriate, we will indicate the
extra results that can be obtained in the proper discrete case.
\nl\nl

\bf 3. Quasi-discrete quantum groups \rm
\nl
Let us start with a pair $(A,\Delta)$ where $A$ is a C$^\ast$-algebra and
$\Delta$ a comultiplication on $A$.  We make the following assumption :

\inspr{3.1}  Assumption \rm
Assume that there is a non-zero element $h$ in $A$ such that
$\Delta(a)(1 \otimes h) = a \otimes h$ for all $a \in A$.
\einspr

If we multiply the above equation to the right with $h^\ast$, we see
immediately that we can assume that $h\geq 0$.
\snl
If we also assume that $\Delta(A)(A \otimes 1)$ is a dense subspace of $A
\overline{\otimes} A$, we can prove a number of results.

\inspr{3.2} Lemma \rm  There exists a non-zero homomorphism $\epsilon : A \rightarrow \Bbb
C$ such that $ah  = \epsilon (a) h$ for all $a \in A$.

\snl\bf Proof : \rm  For all $b \in A$ we have
$$(\Delta(a)(b \otimes 1))(1 \otimes h) = ab \otimes h.$$
By the density of $\Delta(A)(A \otimes 1)$ in $A \overline{\otimes} A$, it
follows that
$$(A \otimes A)(1 \otimes h) \subseteq A \otimes h$$
and so $Ah \subseteq \Bbb Ch$.  This implies the existence of a linear map
$\epsilon : A \rightarrow \Bbb C$ given by $ah = \epsilon (a)h$.  Clearly
$$\epsilon (ab) h = abh = a \epsilon (b) h = \epsilon (b) \epsilon (a) h$$
so that $\epsilon$ is a homomorphism.  Also
$h^\ast h = \epsilon (h^\ast) h$
and because $h \ne 0$ we must have $\epsilon (h^\ast) \ne 0$.
So $\epsilon$ is non-zero.

\inspr{3.3}   Lemma \rm  We can choose $h$ such that $h^2 = h = h^\ast$.

\snl\bf Proof : \rm  We have $h^\ast h = \epsilon (h^\ast)h$.  For any
$\lambda \in \Bbb C$ we
get $(\lambda h)^\ast (\lambda h) = \bar \lambda \epsilon (h^\ast)(\lambda h)$.
So, if $\bar  \lambda \epsilon (h^\ast) = 1$, and if we replace $h$ by $\lambda
h$, we find an element $h$, still satisfying assumption 3.1 but now also $h^\ast
h = h$.  Hence $h = h^\ast$ and $h^2 = h$.
\einspr

From now on, we make this choice of $h$.  In particular $\epsilon(h) = 1$.

\inspr{3.4} Lemma \rm  $\epsilon$ is a $^\ast$-homomorphism.

\snl\bf Proof \rm  For all $a \in A$ we have
$$\align
ha^\ast h &= h \epsilon (a^\ast) h = \epsilon (a^\ast)h\\
ha^\ast h &= (ah)^\ast h = \overline{\epsilon (a)} h\endalign$$
so that $\epsilon (a^\ast) = \overline{\epsilon(a)}.$

\inspr{3.5}  Lemma \rm  For all $a \in A$ we have
$$\align(\iota \otimes \epsilon)\Delta (a) &= a \\
(\epsilon \otimes \iota) \Delta(a) &=a.\endalign$$

\snl\bf Proof : \rm  We have
$$
a \otimes h  = \Delta (a)(1 \otimes h)
 = (\iota \otimes \epsilon)\Delta (a) \otimes h
$$
so that $(\iota \otimes \epsilon)\Delta (a) = a$.
\snl
To prove the second formula, let $a,b \in A$.  Then
$$\split
((\Delta \otimes \iota)(a \otimes b))(1 \otimes h \otimes 1) & = (\Delta (a)
\otimes b)(1 \otimes h \otimes 1) \\
& = (\Delta (a) (1 \otimes h)) \otimes b = a \otimes h \otimes b.
\endsplit
$$
This last element is $\sigma_{23} (a \otimes b \otimes h)$ where $\sigma_{23}$
flips the last two factors.
\snl
If we replace $a \otimes b$ by $\Delta(a)$, we get
$$\split
\sigma_{23} (\Delta (a) \otimes h) & = ((\Delta \otimes \iota)\Delta(a))(1
\otimes h \otimes 1) \\
& = ((\iota \otimes \Delta)\Delta(a)) (1 \otimes h \otimes 1).
\endsplit
$$
Now, we multiply with $c$ in the first factor of the tensor product.  Then we
obtain
$$\sigma_{23} ((\Delta(a)(c \otimes 1)) \otimes h) = (\iota \otimes
\Delta)(\Delta(a)(c \otimes 1))(1 \otimes h \otimes 1).$$
By density of $\Delta(A) (A \otimes 1)$ in $A \overline\otimes A$,
we can replace $\Delta
(a) (c \otimes 1)$ by $a \otimes b$ and we obtain
$$\sigma_{23} (a \otimes b \otimes h) = (a \otimes \Delta (b)) (1 \otimes h
\otimes 1).$$
Hence
$$h \otimes b  = \Delta (b) (h \otimes 1)
 = h \otimes (\epsilon \otimes \iota)\Delta(b)$$
and $b = (\epsilon \otimes \iota)\Delta(b).$
\einspr

Remark that we find that also $\Delta(a)(h \otimes 1) = h \otimes a$ and not
only $\Delta (a) (1 \otimes h) = a \otimes h$.
\snl
From the formulas in 3.5 we also get the uniqueness of $h$.

\inspr{3.6} Lemma \rm  $h$ is unique.

\snl\bf Proof : \rm  Assume that $h$ and $h^\prime$ are non-zero elements satisfying
$\Delta(a)(1 \otimes h) = a \otimes h$ and $\Delta (a) (1 \otimes h^\prime) = a
\otimes h^\prime$.  Assume that $h^2 = h = h^\ast$ and similarly for $h^\prime$.
Let $\epsilon$ and $\epsilon^\prime$ be the associated $^\ast$-homomorphisms.  For
all $a$ we have
$$\epsilon (a)  = \epsilon ((\iota \otimes \epsilon^\prime)\Delta(a))
 = \epsilon^\prime ((\epsilon \otimes \iota)\Delta (a))
 = \epsilon^\prime (a).
$$
So $\epsilon = \epsilon^\prime$.  Then $h = \epsilon (h^\prime)h = h^\prime h =
\epsilon (h) h^\prime = h^\prime$.
\einspr

We summarise in the following proposition.

\inspr{3.7} Proposition \rm  Let $A$ be a C$^\ast$-algebra and $\Delta$ a comultiplication
on $A$.  Assume that $\Delta(A)(A \otimes 1)$ is a dense subspace in $A
\overline
\otimes A$.  Assume there is a non-zero element $h$ in $A$ such that
$\Delta(a)(1 \otimes h) = a \otimes h$ for all $a$.  Then there is a unique
non-zero self-adjoint projection $h$ in $A$ such that $\Delta (a) (1 \otimes h)
= a \otimes h$.  We also have $\Delta (a)(h \otimes 1) = h \otimes a$.  There is
also a unique $^*$-homomorphism $\epsilon : A \rightarrow \Bbb C$ such that $(\iota
\otimes \epsilon)\Delta(a) = a$ and  $(\epsilon \otimes \iota)\Delta(a) = a$
for all $a$.  Furthermore $ah = ha = \epsilon (a) h$ for all $a \in A$.
\einspr

Remark that the statement about the uniqueness of $\epsilon$ was more or less
proved in the proof of lemma 3.6.  Indeed, if $\epsilon$ and $\epsilon^\prime$
are linear maps such that $(\iota \otimes \epsilon^\prime)\Delta(a) = a$ and
$(\epsilon \otimes \iota) \Delta (a) = a$ for all $a$, then $\epsilon =
\epsilon^\prime$.
\snl
We also remark that the condition $\Delta(A)(A \otimes 1)$ dense in $A
\overline
\otimes A$ is not really necessary, a weaker condition like
$$\{ (\omega \otimes \iota) \Delta(a) \mid \omega \in A^\ast, a \in A \}$$
is in $A$ and spans a dense subspace of $A$ would be sufficient to carry out
the above proofs.  Having this condition and the existence of $h$, we see that
we have more, namely
$$\align\{ (\epsilon \otimes \iota)\Delta(a) \mid a \in A \} &= A\\
\{ (\iota \otimes \epsilon)\Delta(a) \mid a \in A \} &= A.\endalign$$
In fact, the condition
$$\{ (\iota \otimes \omega) \Delta (a) \mid \omega \in A^\ast, a \in A \}=A$$
is already a consequence of the existence of $h$ alone.  Indeed, if $a \in A$
and $\varphi \in A^\ast$ such that $\varphi(h) = 1$, then, with $\omega =
\varphi (\cdot h)$ we have
$$(\iota \otimes \omega)\Delta(a)
 = (\iota \otimes \varphi)(\Delta(a)(1 \otimes h)) = a  \varphi (h) = a
$$
\snl
We can illustrate  some of the above results with (easy) examples.

\inspr{3.8} Examples \rm  i) If $A$ is a C$^\ast$-algebra and $\Delta(a) = a \otimes 1$
for all $a$.  Then $\Delta$ is a comultiplication.  Any element $h$ satisfies
$\Delta (a) (1 \otimes h) = a \otimes h$.  For any linear functional $\epsilon$
with $\epsilon (1) = 1$ we have $(\iota \otimes \epsilon)\Delta(a) = a$.  But
here $\Delta(A) (A \otimes 1) = A \otimes 1$ and this is not a dense subspace of
$A \overline \otimes A$.
\snl
ii) If $\Delta(a) = 1 \otimes a$, there is no such
element $h$ while  $\Delta(A)(A
\otimes 1)$ is dense in $A \overline \otimes A$.
\einspr

Now let $A$ be a C$^\ast$-algebra with a comultiplication $\Delta$ and assume
that $\Delta(A)(A \otimes 1)$ and $\Delta(A)(1 \otimes A)$ are dense subspaces
of $A \overline \otimes A$.  We have not yet used the second density
condition, but we will need it later. We know, from our discussion on
the abelian case in section 1 that these two density conditions are
quite natural.
\snl
Assume further that there exists an element $h$ like
before, i.e. a self-adjoint projection $h$ satisfying $\Delta(a)(1 \otimes h) =
a \otimes h$ and $\Delta(a)(h \otimes 1) = h \otimes a$ for all $a$.  Then we
have the existence of a counit $\epsilon$, i.e. a $^\ast$-homomorphism $\epsilon
: A \rightarrow \Bbb C$ such that $(\iota \otimes \epsilon)\Delta(a) = (\epsilon
\otimes \iota)\Delta(a) = a$ for all $a$.  And we have $ah = ha = \epsilon(a)h$.
\snl
In section 1, we also saw that, in the abelian case, we need some extra condition to
distinguish the group case from the case of a semi-group with cancellation.  The
condition that we needed was $\Delta(a)(1 \otimes b) = 0 \Rightarrow a \otimes b
= 0.$
\snl
Here we have the following interesting property.

\inspr{3.9} Proposition \rm  Consider the linear map $T_1 : A \otimes A \rightarrow A \overline
\otimes A$ defined by $T_1(a \otimes b) = \Delta(a)(1 \otimes b)$.  If
$\Delta(h)(1 \otimes a) = 0$ implies $a = 0$ then this map is injective.

\snl\bf Proof : \rm  We know that $(1 \otimes h)\Delta(a) = a \otimes h$.  If we apply $\iota
\otimes \Delta$ and multiply with $1 \otimes 1 \otimes b$, we get
$$(1 \otimes \Delta(h))\Delta^{(2)}(a)(1 \otimes 1 \otimes b) = (a \otimes
\Delta(h))(1 \otimes 1 \otimes b)$$
where we use the notation $\Delta^{(2)}(a)$ for
$(\Delta\otimes\iota)\Delta(a)$ as usual.
This can be rewritten as
$$(1 \otimes \Delta(h))(\Delta \otimes \iota)(T_1(a \otimes b)) = (1 \otimes
\Delta(h))\sigma_{12} (1 \otimes a \otimes b).$$
By linearity we have for all $x \in A \otimes A$,
$$(1 \otimes \Delta(h))(\Delta \otimes \iota)(T_1(x)) =(1 \otimes
\Delta(h))\sigma_{12} (1 \otimes x).$$
If now $T_1(x) = 0$ then $(1 \otimes \Delta(h))\sigma_{12} (1 \otimes x) = 0.$  If
we apply $\omega \otimes \iota \otimes \iota$ with $\omega \in A^\ast$, we find
$\Delta(h)(1 \otimes y) = 0$ with $y = (\omega \otimes \iota)x$.  By assumption
$y = 0$.  Because this is true for all $\omega$, we get also that $x = 0$.
Hence $T_1$ is injective.
\einspr

We see from the above proof that, if $T_1$ can be closed (as a linear map from $A
\overline \otimes A$ to $A \overline \otimes A$), then this closure will still be
injective.  In any case, we see that, if $(x_n)$ is a sequence in $A \otimes A$
such that $T_1x_n \rightarrow 0$ and such that $(x_n)$ converges, then $x_n
\rightarrow 0$.  So, in fact, the inverse of $T_1$ is closable.
\nl
All this leads us to the main object in this paper.

\iinspr{3.10}  Definition \rm  Let $A$ be a C$^\ast$-algebra with a comultiplication
$\Delta$.  Assume that $\Delta(A)(1 \otimes A)$ and $\Delta(A)(A  \otimes 1)$
are dense subspaces of $A \overline \otimes A$.  Assume that $A$ has a projection
$h$ as before and that $\Delta(h)(1 \otimes a) = 0$ implies $a= 0$.
  Then we call $(A,\Delta)$ a
quasi-discrete quantum group.
\einspr

We will see later (in proposition 3.17 below), that it will follow
automatically that also $\Delta(h)(a\ot 1)=0$ implies $a=0$. This in
turn, will give that the map $T_2 : A\ot A \to A\overline\ot A$, defined
as before by $T_2(a\ot b)=\Delta(a)(b\ot 1)$ is injective. The proof is
similar as the one for $T_1$. But it also follows by symmetry. So the
two maps $T_1$ and $T_2$ are injective (compare with [24]). Therefore,
we also get that a quasi-discrete quantum group, where the underlying
C$^*$-algebra is a direct sum of full matrix algebras, is a discrete
quantum group.
\snl
In the remaining of this section, we introduce some subsets of $A$,
 canonically associated to such a quasi-discrete quantum group. We also
 prove some more properties of the element $h$.

\iinspr{3.11} Notation \rm  Let
$$\align
   I_0 &= \{ (\omega \otimes \iota)\Delta(h) \mid \omega \in A^\ast\}\\
   J_0 &= \{ (\iota \otimes \omega)\Delta(h) \mid \omega \in A^\ast \}.
\endalign$$
Let also $I = \overline {I_0}$ and $J = \overline {J_0}$.
\einspr

Remark that any element $\omega\in A^*$ is of the form $\omega_1(\,\cdot\,
a)$ (see e.g. [21]). Hence, by our assumptions, $I_0$ and $J_0$ are
subsets of $A$.
\snl
 We will show that $I=J$ and and that $J$ is a closed two-sided ideal.
We need a few more properties of $h$ before we can do this. These
properties will also be used later.

\iinspr{3.12} Proposition \rm
$\Delta(h)(A\otimes 1)\subseteq\Delta(h)(1\otimes A)^-$

\snl\bf Proof : \rm  Let  $a \in  A$. Because $h^2 = h$ we have,
$$\Delta(h) (a \otimes 1) = \Delta(h) \Delta(h) (a \otimes 1).$$
Since $\Delta (h) (a \otimes 1) \in A \overline \otimes A$ and because $\Delta (A)
(1 \otimes A)$ is dense in $A  \overline \otimes A$, we can approximate $\Delta(h)(a
\otimes 1)$ by finite linear combinations of elements of
the form $\Delta(c)(1 \otimes d)$.
Now
$$
\Delta(h)\Delta(c)(1 \otimes d)  = \Delta (hc)(1 \otimes d)
 = \Delta (h)(1 \otimes \epsilon (c)d). $$
So we see that
$$\Delta(h)(a \otimes 1) = \lim_n \Delta(h) (1 \otimes b_n)$$
for some sequence $(b_n)$ in $A$.

\iinspr{3.13} Proposition \rm  $I$ and $J$ are two-sided ideals.

\snl\bf Proof : \rm
Given $a$ in $A$, we have by the previous result that
$$\Delta(h)(a\otimes 1)=\lim_n\Delta(h)(1\otimes b_n)$$
 for some sequence
$(b_n)$ in $A$.
If we apply $\iota\otimes\omega$ with $\omega \in A^\ast$ we find
$$((\iota \otimes \omega)\Delta(h))a = \lim_n (\iota \otimes \omega_n)\Delta(h)$$
where $\omega_n = \omega (\,\cdot\, b_n)$ for all $n$.  This shows that
$J_0 A
\subseteq J$.  Hence $JA \subseteq J$.  Furthermore $J_0$ and $J$ are clearly
self-adjoint.  A similar argument works for $I$.
\einspr
The injectivity condition gives the following property of the ideal $I$.

\iinspr{3.14} Proposition \rm  The ideal $I$ is essential, i.e. if
$a \in A$ and $ba = 0$ for all $b \in I$, then $a = 0$.

\snl\bf Proof : \rm  If $ba = 0$ for all $b \in I$ then
$$(\omega\ot\iota )(\Delta(h)(1\ot a)) = ((\omega\ot\iota)
\Delta (h)) \cdot a = 0$$
for all $\omega \in A^\ast$.  So $\Delta(h)(1 \ot a) = 0$ and by assumption
$a = 0$.
\einspr

This guarantees that the ideal $I$ carries enough information. In
section 5, where we discuss the Haar measure, we will come back to the
fact that possibly $I$ is strictly smaller than $A$.
\snl
In the discrete case however, the fact that $I$ is essential implies
that it is all of $A$. Moreover, the ideal $\Cal A$ that we defined
before, consisting of the elements $(x_{\alpha})$ with only finitely
many $x_{\alpha}$ non-zero, must be contained in $I_0$. This
can be seen as follows. Take any $\alpha$ and consider the set
$(\omega\ot\iota)\Delta(h)$ with $\omega\in A$ but with support in
$A_{\alpha}$. As in the proof of proposition 3.13, this is also an
ideal. It must be contained in $\Cal A$. Hence, it must be a direct sum
of a finite number of summands of $A$. We must get all components this
way because the set $I_0$ is dense. So we find $\Cal A\subseteq I_0$.
\snl
We now turn back to the general quasi-discrete case. We will show
that $I=J$. To prove this, we
first need some more results on $h$.

\iinspr{3.15} Proposition \rm
 We have
$$\Delta_{23} (h) \Delta_{14} (h) = \Delta_{23} (h) \Delta^{(3)} (h),$$
where we use the `leg numbering' notation (e.g. $\Delta_{23} (h) = 1 \otimes
\Delta (h) \otimes 1$).  We also use $\Delta^{(3)}(h)$ for $(\Delta \otimes \iota
\otimes \iota)(\Delta \otimes \iota)\Delta(h)$ as usual.

\snl\bf Proof : \rm
$$\split
\Delta_{23} (h) \Delta^{(3)} (h) & = (1 \otimes \Delta (h) \otimes 1)(\iota
\otimes \Delta \otimes \iota)(\Delta^{(2)}(h)) \\
& = (\iota \otimes \Delta \otimes \iota)((1 \otimes h \otimes 1)
\Delta^{(2)}(h)).
\endsplit
$$
Now
 $(1 \otimes h \otimes 1)\Delta^{(2)}(a) = (1 \otimes h \otimes
1) \Delta_{13} (a)$.  Therefore
$$\split
\Delta_{23} (h) \Delta^{(3)} (h) & = (\iota \otimes \Delta \otimes \iota)(1
\otimes h \otimes 1)(\Delta_{13} (h)) \\
& = \Delta_{23} (h) \Delta_{14} (h).
\endsplit
$$
\einspr

From this formula we can obtain the following.

\iinspr{3.16} Proposition \rm
We have
$$\Delta_{23} (h) \Delta_{14} (h) \Delta_{12} (h) =
\Delta_{23} (h) \Delta_{14} (h) \Delta_{34} (h).$$

\snl\bf Proof : \rm  We get on the one side
$$\align
   \Delta_{23}(h)\Delta_{14}(h)\Delta_{12}(h)
       &=\Delta_{23}(h)\Delta^{(3)}(h)\Delta_{12}(h)\\
       &=\Delta_{23}(h)(\Delta\ot\iota\ot\iota)(\Delta^{(2)}(h)(h\ot
       1\ot 1))\\
       &=\Delta_{23}(h)(\Delta\ot\iota\ot\iota)(h\ot\Delta(h))\\
       &=\Delta_{23}(h)(\Delta(h)\ot\Delta(h)).
   \endalign$$
On the other side we get
$$\align
   \Delta_{23}(h)\Delta_{14}(h)\Delta_{34}(h)
       &=\Delta_{23}(h)\Delta^{(3)}(h)\Delta_{34}(h)\\
       &=\Delta_{23}(h)(\iota\ot\iota\ot\Delta)(\Delta^{(2)}(h)(1\ot
       1\ot h))\\
       &=\Delta_{23}(h)(\iota\ot\iota\ot\Delta)(\Delta(h)\ot h)\\
       &=\Delta_{23}(h)(\Delta(h)\ot\Delta(h)).
   \endalign$$
This proves the formula.

\iinspr{3.17} Proposition \rm
The two ideals $I$ and $J$ are equal.

\snl\bf Proof : \rm
Consider the formula in the previous proposition. Apply $\iota\ot\iota
\ot\omega\ot\omega_1$ and assume that $\omega|J=0$. Because of the factor
$\Delta_{34}(h)$ in the right hand side of this equation, the result will
be $0$. So we obtain $(a\ot b)\Delta(h)=0$ where $a=(\iota\ot\omega_1)
\Delta(h)$ and $b=(\iota\ot\omega)\Delta(h)$. Because this is true for
all $\omega_1$, we have that this is also true for all $a\in J$.
Hence $(1\ot b)\Delta(h)=0$. By our injectivity assumption, this
implies $b=0$. Because $b=(\iota\ot\omega)\Delta(h)$, we will have that
$\omega|I=0$. This implies that $I\subseteq J$.
\snl
We have shown that the  ideal  $I$ is essential. Because
$I\subseteq J$, this is also the case for $J$. This is equivalent with
the other injectivity property
$$\Delta(h)(a\ot 1)=0 \Rightarrow a=0.$$ Then, by symmetry, or by a
 similar argument, using the factor $\Delta_{12}(h)$ on the left side,
 we get that also $J\subseteq I$.
\einspr

We see above that one injectivity assumption gives the other. This was
also true in the abelian case (see section 1). For the discrete case,
this implies that also $J_0$ is dense in $A$ and in fact, that $\Cal
A\subseteq J_0$.
\snl
It is more or less clear from the definitions that $\Delta(h)(1 \otimes A)
\subseteq J \overline \otimes J$ and that $\Delta(h)(A \otimes 1)
\subseteq J \overline
\otimes J$.  Also $\Delta (h) \in M(J \overline \otimes J)$.
Because $\Delta(h)(1\ot h)=h\ot h$, it follows that $h\in J$.

\iinspr{3.18} Proposition \rm
We have that $\Delta(J)(J\ot 1)$ and $\Delta(J)(1\ot J)$ are dense
subsets of $J\overline\ot J$.

\snl\bf Proof : \rm
We first show that $\Delta(J)(J\ot 1)\subseteq J\overline\ot J$.
For all $\omega$ and $\omega'$ we have
  $$(\jo\ot\omega')\Delta((\jo\ot\omega)\Delta(h))
     =(\jo\ot\omega'\omega)\Delta(h)\in J.$$
Recall that $\omega'\omega$ is defined by
$(\omega'\omega)(a)=(\omega'\ot\omega)\Delta(a)$.
We see that $(\iota\ot\omega')\Delta(J)\subseteq J$ for all $\omega'\in
A^*$. This implies that $\Delta (J)(1\ot J) \subseteq J \overline\ot J$.
Similarly $\Delta (J)(J\ot 1)\subseteq J\overline\ot J$.
\snl
Now we show that $\Delta(J)(J\ot 1)$ is dense in $J\overline\ot J$.
Take $a\in J$. We have $\Delta(a)(1\ot h)=a\ot h$ and hence
$\Delta^{(2)}(a)(1\ot\Delta(h))=a\ot \Delta(h)$.
If we apply $\jo\ot\jo\ot\omega$ we find, with
$b=(\jo\ot\omega)\Delta(h)$, that
   $$a\ot b=(\jo\ot\jo\ot\omega)(\Delta^{(2)}(a)(1\ot\Delta(h))).$$
We claim that this last element can be approximated by elements in
$\Delta(J)(1\ot J)$.  To see this, we work in the universal
representation. Then $\omega=\omega_{\xi,\eta}$ for some vectors $\xi$
and $\eta$ in the underlying Hilbert space. If $(e_i)_i$ is an
orthonormal basis, then
 $$(\jo\ot\jo\ot\omega)(\Delta^{(2)}(a)(1\ot\Delta(h)))
    =\sum_i ((\jo\ot\jo\ot\omega_{e_i,\eta})\Delta^{(2)}(a))
            ((\jo\ot\jo\ot\omega_{\xi,e_i})(1\ot\Delta(h))).$$
This sum converges in norm. Now, these terms in this sum are elements in
$\Delta(J)(1\ot J)$. This gives the density of $\Delta(J)(1 \ot J)$ in
$J\overline\ot J$. The other density is proved in a similar way.
\einspr
\nl\nl

\bf 4. The antipode in quasi-discrete quantum groups
\rm
\bigskip
In this section, the pair $(A,\Delta)$ is a quasi-discrete quantum group and $h$
denotes the unique projection such that $\Delta(a)(1 \otimes h) = a \otimes h$
for all $h$.  We know that also $\Delta(a)(h \otimes 1) = h \otimes a$ for all
$a$.  By the injectivity assumptions $\Delta(h)(1 \otimes a) = 0$ implies $a =0$
and $\Delta(h)(a \otimes 1) = 0$ implies $a = 0$, we know that, in a way,
$\Delta(h)$ carries enough information.
\snl
We now want  to define the antipode.  In the case of a Hopf algebra, we have
$$a \otimes 1 = \sum_{(a)} a_{(1)} \otimes a_{(2)} S(a_{(3)})$$
where we use the common notational convention.  So, if such an element $h$ exists,
we get
$$\split
\Delta(h)(a \otimes 1) & = \sum_{(a)} \Delta (h) \Delta(a_{(1)})(1  \otimes
S(a_{(2)})) \\
& = \sum_{(a)} \Delta (h a_{(1)}) (1 \otimes S(a_{(2)})) \\
& = \sum_{(a)} \epsilon (a_{(1)}) \Delta(h)(1 \otimes S(a_{(2)})) \\
& = \Delta (h) (1 \otimes S(a)).
\endsplit
$$
We will use this formula to define $S$.  Recall that $\Delta(h)(a \otimes 1)
\subseteq (\Delta(h)(1 \otimes A))^-$ when we have a quasi-discrete quantum
group (Proposition 3.12).

\inspr{4.1} Definition \rm  Define the antipode $S$ on $A$ by
$$\Cal D (S) = \{ a \in A \mid \exists b \in A \ \text  { such that }\ \Delta(h) (a
\otimes 1) = \Delta(h)(1 \otimes b) \}$$
and $S(a)$ by $\Delta(h)(a \otimes 1) = \Delta (h)(1 \otimes S(a)).$
\einspr

Remark that the element $b$ above is unique because of our injectivity
assumptions.  So $S$ is well defined.  By the same argument $S$ is injective.

\nl
Because $\Delta(h)(h\ot 1)=h\ot h=\Delta(h)(1\ot h)$, it
is immediately clear from the definition that $h\in\Cal D(S)$ and
that $S(h)=h$. On the other hand, it is not clear if there are enough
elements in $\Cal D(S)$. We will show later that this is the case. Now
we prove some properties of $S$ that follow easily from the definition.

\inspr{4.2} Proposition \rm  $S$ is closed.

\bf\snl Proof : \rm  Let $(a_n)$ be a sequence in $\Cal D(S)$ and assume that $a_n
\rightarrow a$ and $S(a_n) \rightarrow b$ with $a, b\in A$.  Then $\Delta(h)(a_n \otimes 1) =
\Delta(h)(1 \otimes S(a_n))$ for all $n$ and in the limit we get $\Delta(h)(a
\otimes 1) = \Delta(h)(1 \otimes b)$.  This shows that $a \in \Cal D(S)$ and
that $b = S(a)$.  So $S$ is closed.

\inspr{4.3} Proposition \rm  If $a,b \in \Cal D(S)$, then $ab \in \Cal D(S)$ and $S(ab) =
S(b)S(a)$.

\snl\bf Proof : \rm  If $a,b \in \Cal D(S)$ we have
$$\split
\Delta(h)(ab \otimes 1) & = \Delta(h)(a \otimes 1)(b \otimes 1) \\
& = \Delta(h)(1 \otimes S(a))(b \otimes 1) \\
& = \Delta(h)(b \otimes 1)(1 \otimes S(a)) \\
& = \Delta(h)(1 \otimes S(b) S(a)).
\endsplit
$$
This shows that $ab \in \Cal D(S)$ and that $S(ab) = S(b)S(a).$
\einspr

We want to show that $S(a)^\ast \in \Cal D(S)$ when $a \in \Cal D(S)$ and that
$S(S(a)^\ast)^\ast = a$.  Before we can do this, we need two more results.

\inspr{4.4} Proposition \rm  For all $a \in A$ we have
$$(\Delta(h) \otimes 1)(a \otimes 1 \otimes 1)(1 \otimes \Delta (h)) =
(\Delta(h) \otimes 1)(1 \otimes 1 \otimes a)(1 \otimes \Delta(h)).$$

\snl\bf Proof : \rm  Given $a \in A$ we have
$$\split
(\Delta (h) \otimes 1)  (a \otimes 1 \otimes 1)(1 \otimes \Delta(h))
& = (\Delta (h) \otimes 1)(\iota \otimes \Delta) (a \otimes h) \\
& = (\Delta (h) \otimes 1)(\iota \otimes \Delta) (\Delta(a)(1 \otimes h)) \\
& = (\Delta (h) \otimes 1) (\Delta \otimes \iota)(\Delta (a))(1 \otimes
\Delta(h)) \\
& = (\Delta \otimes \iota)((h \otimes 1)\Delta(a))(1 \otimes \Delta(h)) \\
& = (\Delta \otimes \iota)(h \otimes a)(1 \otimes \Delta(h)) \\
& = (\Delta (h) \otimes 1)(1 \otimes 1 \otimes a)(1 \otimes \Delta(h)).
\endsplit
$$

\inspr{4.5} Proposition \rm  If $a \in \Cal D(S)$ then also
$$(1 \otimes a)\Delta(h) = (S(a) \otimes 1)\Delta(h).$$

\snl\bf Proof : \rm  By the previous proposition  we know that
$$(\Delta(h) \otimes 1)(1 \otimes S(a) \otimes 1)(1 \otimes \Delta(h)) =
(\Delta(h) \otimes 1)(1 \otimes 1 \otimes a)(1 \otimes \Delta(h)).$$
If we apply $\iota \otimes \iota \otimes \omega$ with $\omega \in A^\ast$ we get
$$\Delta(h)(1 \otimes x) = \Delta(h)(1 \otimes y)$$
with $x = (\iota \otimes \omega)(S(a) \otimes 1)\Delta(h))$ and
$y = (\iota \otimes \omega)((1 \otimes a)\Delta(h))$.
By the injectivity assumption, we have $x = y$.  And since this holds for all
$\omega$, we get
$$(S(a) \otimes 1) \Delta (h) = (1 \otimes a)\Delta(h).$$
\einspr

We can verify this formula  in the case of a Hopf algebra.
Indeed
$$\split
\sum_{(a)} (S(a_{(1)})a_{(2)} \otimes a_{(3)}) \Delta (h) & = \sum_{(a)}
(S(a_{(1)}) \otimes 1)\Delta(a_{(2)} h) \\
& = \sum_{(a)} (\epsilon(a_{(2)})S(a_{(1)}) \otimes 1)\Delta(h) \\
& = (S(a) \otimes 1)\Delta(h).
\endsplit
$$
On the other hand, this expression is equal to
$$\sum_{(a)} (\epsilon(a_{(1)}) 1 \otimes a_{(2)}) \Delta (h) = (1 \otimes a)
\Delta (h).$$
\snl
From the formula in 4.5 we can easily proof the
formula $S(S(a)^\ast)^\ast = a$.

\inspr{4.6}  Proposition \rm  If $a \in \Cal D(S)$ then $S(a)^\ast \in \Cal D(S)$ and
$S(S(a)^\ast)^\ast = a$.

\snl\bf Proof : \rm  For $a \in \Cal D(S)$ we have
$$(1 \otimes a)\Delta(h) = (S(a) \otimes 1)\Delta(h).$$
If we take adjoints, we obtain
$$\Delta(h)(S(a)^\ast \otimes 1) = \Delta(h)(1 \otimes a^\ast).$$
This shows that $S(a)^\ast \in \Cal D(S)$ and that $S(S(a)^\ast) =
a^\ast$.
\einspr

Our next objective is to try to prove the (equivalent of) the formula
$$\Delta(S(a)) = \sigma(S \otimes S)\Delta(a),$$ where $\sigma$ is the flip.
This will follow from the formula in proposition 3.15.

\inspr{4.7} Proposition \rm  If $a \in \Cal D(S)$, then
$$\Delta_{23}(h)\Delta_{14}(h)(\Delta (a) \otimes 1 \otimes 1) = \Delta_{23} (h)
\Delta_{14} (h) (1 \otimes 1 \otimes \Delta (S(a)).
$$

\snl\bf Proof : \rm  If $a \in \Cal D(S)$, then
$$\split
\Delta_{23}(h)\Delta_{14} (h) (\Delta(a) \otimes 1 \otimes 1) & = \Delta_{23}
(h) \Delta^{(3)} (h) (\Delta \otimes \Delta)(a \otimes 1) \\
& = \Delta_{23} (h) (\Delta \otimes \Delta)(\Delta(h)(a \otimes 1)) \\
& = \Delta_{23} (h) (\Delta \otimes \Delta)(\Delta(h)(1 \otimes S(a))) \\
& = \Delta_{23} (h) \Delta^{(3)} (h) (1 \otimes 1 \otimes \Delta (S(a))) \\
& = \Delta_{23} (h) \Delta_{14} (h) (1 \otimes 1 \otimes \Delta(S(a))).
\endsplit
$$
\einspr

We can rewrite this formula as
$$(\Delta(h) \otimes \Delta(h)) (\Delta_{13}(a)) = (\Delta(h) \otimes
\Delta(h))\Delta_{42} (S(a))$$
using the right permutation.  Hence we see that this formula means $(S \otimes
S)\Delta (a) = \sigma \Delta (S(a))$.
\snl
Now, recall the definitions of $I_0$ and $J_0$ (see 3.11).  We had
$$\align I_0 &=\{ (\omega \otimes \iota) \Delta (h) \mid \omega \in A^\ast\}\\
J_0 &= \{ (\iota \otimes \omega) \Delta (h) \mid \omega \in A^\ast \}.
\endalign$$
We also defined the closures $I = \overline{I_0}$ and $J =
\overline{J_0}$
and we saw that $I=J$ and that $I$ is an essential ideal of $A$.
\snl
We will now prove some properties of $\Cal D(S)$ in connection with
$I_0$ and $J_0$.

\inspr{4.8} Proposition \rm   Let $\omega \in A^\ast$ and assume that $\omega(a) = \psi
(\Delta(h)(a \otimes 1))$ for some $\psi \in (A \overline \otimes A)^\ast$.
Then $(\iota
\otimes \omega)\Delta(h) \in \Cal D (S) \cap J_0$ and $S((\iota \otimes
\omega)\Delta (h)) = (\omega_1 \otimes \iota)\Delta(h)$ where $\omega_1 (a) =
\psi (\Delta(h)(1 \otimes a))$.

\snl\bf Proof : \rm  Let $\omega, \omega_1$ and $\psi$ be as in the
formulation of the proposition.
Let $a = (\iota \otimes \omega)\Delta(h)$ and $b = (\omega_1 \otimes
\iota)\Delta(h)$.  Then
$$\split
\Delta (h) (a \otimes 1) & = \Delta (h) ((\iota \otimes \omega) \Delta (h)
\otimes 1) \\
& = (\iota \otimes \omega \otimes \iota) (\Delta_{13} (h) \Delta_{12} (h)) \\
& = (\iota \otimes \psi \otimes \iota) (\Delta_{23} (h) \Delta_{14} (h)
\Delta_{12} (h)).
\endsplit
$$
On the other hand, we get
$$\split
\Delta(h)(1 \otimes b) & = \Delta (h) (1 \otimes (\omega_1 \otimes
\iota)\Delta(h)) \\
& = (\iota \otimes \omega_1 \otimes \iota)(\Delta_{13} (h) \Delta_{23} (h)) \\
& = (\iota \otimes \psi \otimes \iota)(\Delta_{23} (h) \Delta_{14} (h)
\Delta_{34} (h)).
\endsplit
$$
Then, it follows from the formula in proposition 3.16, that these two
expressions are the same.
Therefore $\Delta(h)(a \otimes 1) = \Delta(h)(1 \otimes b)$ and we obtain
$a \in \Cal D(S)$ and $S(a) = b$.
\einspr

In the next proposition, we will see that there are, in a way, enough
elements in $\Cal D(S)$.

\inspr{4.9}  Proposition \rm  $\Cal D(S) \cap J_0$ is dense in $J_0$.

\snl\bf Proof : \rm  We saw that $(\iota \otimes \omega)\Delta (h) \in \Cal D(S) \cap J_0$ if
$\omega(a) = \psi (\Delta (h)(a \otimes 1))$ for some $\psi \in
(A \overline{\otimes} A)^\ast$.  Now assume that $\varphi \in A^\ast$ and that $\varphi((\iota \otimes
\omega) \Delta (h)) = 0$ for all such $\omega$.  Then $\omega(a) = 0$ for all
such $\omega$ when $a = (\varphi \otimes \iota)\Delta(h)$.  This means that
$\psi (\Delta(h)(a \otimes 1))=0$ for all $\psi \in (A \overline{\otimes}
 A)^\ast $.  Hence
$\Delta(h)(a \otimes 1)=0$ and by the injectivity assumption, $a = 0$.
So we find
$\varphi ((\iota \otimes \omega)\Delta(h)) = \omega (a) = 0$ for all $\omega \in
A^\ast$.  This means that $\varphi | J_0 = 0$.  So $\Cal D(S) \cap J_0$ is dense
in $J_0$.
\einspr

It is not clear what happens outside $J_0$.  If we look at the elements we get
in $\Cal D(S) \cap J_0$ and if we look at the image under $S$, we find elements
in $I_0$.  If we take the adjoints, we find elements in $\Cal D(S) \cap I_0$.
In fact, we can also show, in a similar way, that $\Cal D(S) \cap I_0$ is dense
in $I_0$.
\nl
We can also prove the following.

\iinspr{4.10} Proposition \rm  $J_0 \Cal D(S) \subseteq J_0$ and $\Cal D(S) I_0 \subseteq
I_0$.

\snl\bf Proof : \rm  Because $\Delta(h)(a \otimes 1) = \Delta(h)(1 \otimes S(a))$ when $a \in
\Cal D(S)$ we get
$$((\iota \otimes \omega)\Delta(h))a = (\iota \otimes \omega_1)(\Delta(h))$$
where $\omega_1(x) = \omega (xS(a))$.  This proves that $J_0 \Cal D(S) \subseteq
J_0$.  The other result is proved by using $(1 \otimes a)\Delta(h) = (S(a)
\otimes 1)\Delta(h)$.

\iinspr{4.11} Proposition \rm  $\Cal D(S) \cap J_0$ is a subalgebra of $\Cal D(S)$.

\snl\bf Proof : \rm
$$\align(\Cal D(S) \cap J_0)(\Cal D(S) \cap J_0)
& \subseteq \Cal D(S) \Cal D(S) \subseteq \Cal D(S)\\
(\Cal D(S) \cap J_0)(\Cal D(S) \cap J_0)
& \subseteq J_0\Cal D(S) \subseteq J_0.\endalign$$
\einspr

We will need some more results of this type when we treat the regular
representation and the Haar measure.
\snl
Let us now discuss the discrete case. We have seen that $J_0$ contains
all summands of $A$. It follows from the previous results that $\Cal A$
is a subalgebra of $\Cal D(S)$. By taking adjoints, it is also a
subalgebra of $\Cal D(S^{-1})$. In fact, from the proof of 4.8, we see
that $S(\Cal A)=\Cal A$. This is quite normal. The pair $(\Cal
A,\Delta)$ is a multiplier Hopf $^*$-algebra, the antipode exists and
maps $\Cal A$ onto $\Cal A$. Moreover, the antipode is unique and it
must coincide with the antipode that we obtain here.
\snl
In this case, we also have the following phenomenon (see also [6]). If
$e_{\alpha}$ is the identity in $A_{\alpha}$, then $e_{\alpha} \in \Cal
D(S)$ and $S(e_{\alpha})$ is again a minimal central projection. Hence,
it must be some $e_{\alpha'}$. By the fact that
$S(S(e_{\alpha})^*)^*=e_{\alpha}$ we find $S(e_{\alpha'})=e_{\alpha}$. So,
although $S^2\neq \iota$ may occur, we do have that $S$ is involutive on
the indices and on the components of $A$. The dimensions of $A_{\alpha}$
and $A_{\alpha'}$ must be the same. We see further that
$A_{\alpha'}=\{ (\iota\ot\omega)\Delta(h) \mid \omega\in A^*$ with
support in $A_{\alpha} \}$. See also the remark after propositiong 3.14.
\nl
We now want to look at the adjoint of $S$ as an operator on $A^\ast$.  We have to
be carefull since the domain need not be dense.  Formally we must have $(S_0
\omega)(a) = \omega(S(a))$ when $\omega \in \Cal D(S_0)$, when we use $S_0$ for
the operator on $A^\ast$.  If $\omega$ has the form $\omega_1$ as in 4.8., that
is, if
$\omega(a) = \psi (\Delta(h)(1 \otimes a))$
then $(S_0 \omega)(a) = \psi (\Delta (h)(1 \otimes S(a)) = \psi (\Delta (h) (a
\otimes 1)).$  This suggests the following definition.

\iinspr{4.12} Definition \rm  Define a map $S_0$ on $A^\ast$ by
$$\Cal D(S_0) = \{ \omega \in A^\ast \mid \exists \psi \in
(A \overline{\otimes} A)^\ast \text
{ such that } \omega(a) = \psi (\Delta (h) (1 \otimes a)) \}$$
and $(S_0 \omega)(a) = \psi (\Delta (h)(a \otimes 1))$.
\einspr

Remark that $S_0$ is well-defined.  If $\psi (\Delta(h)(1 \otimes a)) = 0$ for
all $a$, then also $\psi (\Delta (h)(a \otimes 1)) = 0$ for all $a$ because, as
we saw before
$$\Delta(h)(a \otimes 1)  \subseteq (\Delta (h) (1 \otimes A))^-.$$
A similar argument will give here that $S_0$ is injective.

\iinspr{4.13} Proposition \rm  $\Cal D(S_0)$ is $w^\ast$-dense in $A^\ast$.

\snl\bf Proof : \rm  If $a \in A$ and $\omega (a) = 0$
for all $\omega \in \Cal D(S_0)$, then
$\psi (\Delta (h)(1 \otimes a)) = 0$ for all $\psi \in (A \overline \otimes A)^\ast$.
Hence $\Delta (h) (1 \otimes a) = 0$ and $a = 0$.
\einspr

From the motivation before the definition, we saw already that the operators $S$
on $A$ and $S_0$ on $A^\ast$ are adjoint to each other in the sense that $(S_0
\omega)(a) = \omega (S(a))$ when $a \in \Cal D(S)$ and $\omega \in \Cal D(S_0)$.
The fact that $\Cal D(S)$ is not dense in $A$ is related with the fact that
$S_0$ need not be closed (or closable) in $A^\ast$.
\snl
We could however restrict $S$ to $\Cal D(S) \cap J_0$ and obtain a linear map
from $J$ to $J$ which is densily defined.  Then, the adjoint would become a map
from $J^\ast$ to $J^\ast$.  We might  loose some information because
 $\Cal D(S)$ could be larger. On the other hand, we will only work with
 this restriction of $S$.
\nl
We finish this section by showing that also $S_0$ is an anti-homomorphism.

\iinspr{4.14} Proposition \rm  If $\omega_1, \omega_2 \in \Cal D(S_0)$, then
$\omega_1\omega_2$, defined by $(\omega_1 \omega_2)(x) = (\omega_1 \otimes
\omega_2) \Delta (a)$, is also in $\Cal D(S_0)$ and $S_0 (\omega_1 \omega_2) =
(S_0(\omega_2))(S_0(\omega_1))$.

\snl\bf Proof : \rm  Let $\psi_1, \psi_2 \in (A \overline\otimes A)^\ast$ be given such that
$$\align\omega_1 (a) &= \psi_1 (\Delta (h)(1 \otimes a))\\
\omega_2 (a) &= \psi_2 (\Delta (h)(1 \otimes a))\endalign$$
for all $a \in A$.
Then
$$(\omega_1 \omega_2)(a) = (\omega_1 \otimes \omega_2)(\Delta(a)) = (\psi_1
\otimes \psi_2)((\Delta (h) \otimes \Delta(h))(\Delta_{24} (a))).$$
On the other hand
$$\split
((S_0 \omega_2)(S_0 \omega_1))(a) & = (S_0 \omega_2 \otimes S_0
\omega_1)(\Delta(a)) \\
& = (\psi_2 \otimes \psi_1)((\Delta(h) \otimes \Delta (h))(\Delta_{13} (a))) \\
& = (\psi_1 \otimes \psi_2)((\Delta(h) \otimes \Delta (h))(\Delta_{31}
(a)).
\endsplit
$$
We have the formula
$$\Delta_{23}(h) \Delta_{14}(h) = \Delta_{23} (h) \Delta^{(3)} (h).$$
If we use the permutation $\sigma$ given by
$$\sigma(a \otimes b \otimes c \otimes d) = b \otimes c \otimes a \otimes d,$$
we see that
$$\split
(\Delta(h) \otimes \Delta(h)) \Delta_{24} (a) & = \sigma (\Delta_{23} (h)
\Delta_{14} (h) \Delta_{34} (a)) \\
& = \sigma (\Delta_{23} (h) \Delta^{(3)} (h) \Delta_{34} (a)) \\
& = (\Delta (h) \otimes 1) \sigma ((\Delta \otimes \Delta)(\Delta (h) (1
\otimes a))).
\endsplit
$$
So $(\omega_1 \omega_2) (a) = \psi (\Delta (h)(1 \otimes a))$
when
$$\psi (x) = (\psi_1 \otimes \psi_2)(\Delta (h) \otimes 1) \sigma ((\Delta
\otimes \Delta)(x)))$$
for $x \in A \overline\otimes A$.  This shows already that $\omega_1 \omega_2 \in
\Cal D(S_0)$.
Similarly,
$$\split
(\Delta (h) \otimes \Delta (h))(\Delta_{31} (a)) & = \sigma (\Delta_{23} (h)
\Delta_{14} (h) \Delta_{12} (a)) \\
& = \sigma (\Delta_{23} (h) \Delta^{(3)} (h) \Delta_{12} (a)) \\
& = (\Delta (h) \otimes 1)\sigma(
(\Delta \otimes \Delta) (\Delta (h) (a \otimes 1)))
\endsplit
$$
so that
$$(S_0 (\omega_2))(S_0 (\omega_1))(a) = \psi (\Delta (h) (a \otimes 1)).$$
This proves that also $S_0(\omega_1 \omega_2) = (S_0 \omega_2)(S_0 \omega_1).$
\einspr

\nl

\bf 5. The Haar measure \rm
\nl
If $(A, \Delta)$ is a compact quantum group, the (right invariant) Haar measure
is a positive linear functional $\varphi$ on $A$ such that $(\varphi \otimes
\omega)\Delta(a) = \varphi(a)\omega(1)$ whenever $\omega \in A^\ast$.  We would
like to use this as a motivation for the definition of a Haar measure on a
quasi-discrete quantum group.
\snl
So, let $(A,\Delta)$ be a quasi-discrete quantum group. In this section
we will assume that $A$ is separable for technical convenience. Now the
Haar measure will be a weight $\varphi$ on
$A$.  Formally, we  need
$$\varphi((\iota \otimes \omega)\Delta(h)) = \varphi(h)\omega(1),$$
or if we normalise $\varphi$ so that $\varphi(h) = 1$, that
$$\varphi((\iota \otimes \omega)\Delta (h))  = \omega(1)$$
for all $\omega \in A^\ast$.  Then, if $a \in A$ and if $a$ has the form $(\iota
\otimes \omega_1)\Delta(h)$ where $\omega_1\in A^\ast$, we will get, again formally,
$$\split
\varphi ((\iota \otimes \omega)\Delta(a)) & = \varphi ((\iota \otimes \omega
\otimes \omega_1)(\Delta^{(2)} (h))) \\
& = \varphi ((\iota \otimes \omega \omega_1)(\Delta(h))) \\
& = (\omega \omega_1)(1) = \omega(1)\omega_1(1) \\
& = \varphi ((\iota \otimes \omega_1)\Delta(h)) \cdot \omega (1) \\
& = \varphi (a) \omega (1).
\endsplit
$$
Apart from the
 technical problems, there are some more fundamental ones here, as we
will point out below.
\snl
The first problem is a consequence of the fact that the ideal $J$ may
not be all of $A$.  So, the element $\omega$  is not determined by the element
$(\iota \otimes \omega)\Delta(h)$.  However, we have the following.

\inspr{5.1} Lemma \rm
  If $\omega, \omega^\prime \in A^\ast$ and $(\iota \otimes
\omega)\Delta(h) = (\iota \otimes \omega^\prime) \Delta(h)$, then $\omega | J =
\omega^\prime | J$.

\snl\bf Proof : \rm
  If we apply any $\omega^{\prime \prime} \in A^\ast$ we get $\omega(a) =
\omega^\prime(a)$ for $a = (\omega^{\prime \prime} \otimes \iota)\Delta(h)$.
Such elements are dense in $J$ and by continuity we have $\omega | J =
\omega^\prime | J$.
\einspr
Then we can define a weight $\varphi$ on $A$.

\inspr{5.2} Definition \rm
  Define $\varphi : A^+ \rightarrow [0,\infty]$ by $\varphi(x) =
\| \omega | J \|$ when $x = (\iota \otimes \omega)\Delta(h)$ and $\omega \in
A^\ast_+$ and let $\varphi(x) = \infty$ when $x \in A^+$ but not of this form.
\einspr

We will prove that this is a weight.  The main problem is to show that $0
\leq x \leq y$ and $y = (\iota \otimes \omega)\Delta(h)$ with $\omega \in
A^\ast_+$ implies that also $x = (\iota \otimes \omega^{\prime})\Delta(h)$ for some
$\omega^\prime \in A^\ast_+$.  This is necessary (and in fact sufficient) to show
that $\varphi(x + y) = \varphi (x) + \varphi (y)$.
\snl
Before we can prove this, we need some other results.  These results will also
be useful in the next section.

\inspr{5.3} Proposition \rm  $J^+_0$ is dense in $J^+$.

\snl\bf Proof : \rm
  We have seen that $\Cal D(S) \cap J_0$ is dense in $J_0$ (proposition
4.9).  We also know that $J_0 \Cal D(S) \subseteq J_0$ (proposition
4.10).
Take any $x \in J^+$.  Consider a sequence $(a_n)$ of elements in $\Cal D(S)
\cap J_0$ such that $a_n \rightarrow x^{1/2}$.  Then $a_n^\ast a_n \rightarrow
x$ and
$$\split
a^\ast_n a_n & \subseteq (\Cal D(S) \cap J_0)^\ast (\Cal D(S) \cap J_0) \\
& \subseteq J_0 \Cal D(S) \subseteq J_0 .
\endsplit
$$
So $x$ is the limit of elements in $J^+_0$.
\einspr

Similarly, $I_0^+$ is dense in $J^+$.

\inspr{5.4} Proposition \rm
  If $\omega \in A^\ast$ and $(\iota \otimes \omega)\Delta(h) \geq
0$, then $\omega | J \geq 0$.

\snl\bf Proof : \rm
Let $\overline \omega (x) = \omega (x^\ast)^-$.  Then
$$(\iota \otimes \overline \omega)\Delta(h) = ((\iota \otimes \omega)\Delta(h))^\ast
= (\iota \otimes \omega) \Delta (h).$$
So, if we replace $\omega$ by $\frac 12 (\omega + \overline \omega)$, we may assume that
$\omega$ is self-adjoint.  Then decompose $\omega = \omega^+ - \omega^-$ where
$\omega^+$ and $\omega^-$ are positive.  We have $\omega \leq \omega^+$ so that
$$0 \leq (\iota \otimes \omega) \Delta (h) \leq (\iota \otimes \omega^+) \Delta
(h).$$
Write $x = (\iota \otimes \omega)\Delta(h)$ and $y = (\iota \otimes
\omega^+)\Delta (h)$.  Put $u_n = {1 \over (y + {1 \over n})^{1/2}} \cdot
x^{1/2}$.  Then
$$u_n u_n^\ast = {1 \over (y + {1 \over n})^{1/2}} x {1 \over (y + {1 \over
n})^{1/2}} \leq {y \over y + {1 \over n}} \leq 1$$
and
$$u^\ast_n y u_n = x^{1/2} {y \over y + {1 \over n}} x^{1/2} \rightarrow x$$
because
$$\split
\| x - u^\ast_n y u_n \| & = \| x^{1/2} (1 - {y \over y + {1 \over n}} )x^{1/2}
\| \\
& = {1 \over n} \| x^{1/2} {1 \over y + {1 \over n}} x^{1/2} \| \\
& = {1 \over n} \| u^\ast_n u_n \| = {1 \over n} \| u_n u^\ast_n \| \leq {1 \over
n}.
\endsplit
$$
We have $x \in J^+_0$.  Hence $x^{1/2} \in J$ and $u_n \in J$.  Choose $v_n \in
\Cal D(S) \cap J_0$ such that $\| v_n - u_n \| \rightarrow 0$.  Then also
$v^\ast_n y v_n \rightarrow x$.  Now we have
$$\split
v^\ast_n y v_n & = v^\ast_n ((\iota \otimes w^+)\Delta(h))v_n \\
& = (\iota \otimes \omega^+)((v^\ast_n \otimes 1)\Delta (h) (v_n \otimes 1)) \\
& = (\iota \otimes \omega^+)((1 \otimes S(v_n)^\ast) \Delta(h) (1 \otimes
S(v_n)) \\
& = (\iota \otimes \omega_n) (\Delta(h))
\endsplit
$$
where $\omega_n(a) = \omega^+ (S(v_n)^\ast aS (v_n))$.
\snl
If we apply $\omega^\prime\in A^\ast$ we get
$$\omega^\prime (x) = \omega ((\omega^\prime \otimes \iota) \Delta (h))$$
on the one hand and
$$\omega^\prime (x) = \lim_n \omega^\prime (v^\ast_n y v_n) = \lim_n\omega_n
((\omega^\prime \otimes 1)\Delta(h)).$$
Therefore $\omega (a) = \lim \omega_n(a)$ for all $a \in I_0$.  If $a \geq 0$
then $\omega_n(a) \geq 0$ so that $\omega (a) \geq 0$.  Because $I^+_0$ is dense
in $J^+$ we get $\omega|J \geq 0$.

\inspr{5.5} Remark \rm
  If $\omega \in A^\ast$ and $(\iota \otimes \omega) \Delta (h) \geq
0$ then, for all $\omega^\prime \in A^\ast$ with $\omega^\prime \geq 0$ we get
$\omega((\omega^\prime \otimes \iota) \Delta (h)) \geq 0$.  We know already that
$I^+_0$ is dense in $J^+$.  If we also knew that $\{ (\omega^\prime \otimes
\iota) \Delta (h) \mid \omega^\prime \geq 0 \}$ is dense in $J^+$, we would obtain
immediately that $\omega | J \geq 0$.  But for the density of $\{ (\omega^\prime
\otimes \iota) \Delta (h) \mid \omega^\prime \geq 0 \}$ in $J^+$, we precisely need
the density of $I^+_0$ in $J^+$ and (the analogue for $I_0$ of) the previous
proposition.
\einspr

So, if we combine 5.3 and 5.4 we get that the set $\{ ( \iota \otimes
\omega)(\Delta (h)) \mid \omega \in A^\ast_+ \}$ is dense in $J_+$.  We need some
more results before we can prove the additive property of the Haar weight.

\inspr{5.6} Proposition \rm
  If $a \in \Cal D(S) \cap J_0$ and $x \in J_+$ then $a^\ast x
a \in J_0$.

\snl\bf Proof : \rm
We can approximate $x^{1/2}$ by a sequence $(b_n)$ in $\Cal D(S) \cap
J_0$ such that $\| b_n \| \leq \| x^{1/2} \|$.  Let $y_n = a^\ast b^\ast_n b_n
a$ for all $n$.  Then $0 \leq y_n \leq \| x \| a^\ast a$.  Because $\Cal D(S)
\cap J_0$ is a subalgebra (proposition 4.11) we have $b_n a \in \Cal D(S) \cap
J_0$.  And, as in 5.3, we know that $(\Cal D(S) \cap J_0)^\ast (\Cal D(S)  \cap
J_0) \subseteq J_0$.  So $y_n \in J^+_0$ and $a^\ast a \in J^+_0$.  We also have
of course that $y_n \rightarrow a^\ast x a$.
\snl
Now let $\omega_n, \omega^\prime \in A^\ast_+$ such that $y_n  = (\iota \otimes
\omega_n) \Delta (h)$ and $\| x \| a^\ast a = (\iota \otimes
\omega^\prime)\Delta(h)$.  Then $0 \leq y_n \leq \| x \| a^\ast a$ will imply
$\omega_n \leq \omega^\prime$ on $J$ by 5.4.  By compactness, we may assume that
 $\omega_n \rightarrow  \omega$ for some $\omega \in A^\ast_+$ in
the $w^\ast$-topology on $J$.  So, for all $\psi \in A^\ast$ we get
$$\split
\psi (y_n) & = (\psi \otimes \omega_n) \Delta (h) \\
& = \omega_n ((\psi \otimes \iota)\Delta (h)).
\endsplit
$$
This converges to
$$\omega((\psi \otimes \iota)\Delta(h)) = \psi ((\iota \otimes
\omega)\Delta(h)).$$
On the other hand $\psi(y_n) \rightarrow \psi(a^\ast xa)$.  It follows that
$a^\ast xa = (\iota
\otimes \omega)\Delta(h)$.
\einspr

Also this result would of course follow easily if we had already the hereditary
property of $J^+_0$. We can prove this now.

\inspr{5.7} Proposition \rm
  If $0 \leq x \leq y$ and $y \in J^+_0$, then $x \in J^+_0$.

\snl\bf Proof : \rm
  Choose $\omega \in A^\ast_+$ such that $y = (\iota \otimes \omega)\Delta
(h)$.  For any $a \in \Cal D(S) \cap J_0$ we have $a^\ast y a = (\iota \otimes
\omega_a)\Delta (h)$ where $\omega_a = \omega (S(a)^\ast \cdot S(a))$.  We also
have $a^\ast xa$ in $J^+_0$.
\snl
Now, take $a$ of the form $a=(\jo\ot\omega)\Delta(h)$ where
$\omega(z)=\psi(\Delta(h)(z\ot 1))$ for some $\psi\in (A\overline\ot
A)^\ast$ and all $z$. Then $S(a)=(\omega_1\ot\jo)\Delta(h)$ where
$\omega_1(z)=\psi(\Delta(h)(1\ot z))$. In the proof of proposition 4.9,
we saw that such elements are dense in $J$. Similarly, the elements
$S(a)$ with $a$ of this form are dense in $J$.  If we consider an
approximate identity in $J$ and approximate these in turn by elements
$S(a)$ as above, we find a sequence $(a_n)$ of elements in $\Cal
D(S)\cap J_0$ such that $S(a_n)\to 1$ in the strict topology of $M(J)$.
We can also assume that $\|S(a_n)\|\leq 1$ for all $n$.
\snl
Then, $\omega_{a_n}$ is a bounded sequence and we may assume that it
converges to $\omega'$. Then
   $$a_n^\ast y a_n
   =(\jo\ot\omega_{a_n})\Delta(h)\to(\jo\ot\omega')\Delta(h)$$
weakly. On the other hand
$$\split
\Delta_{13} (h)(1 \otimes 1 \otimes a^\ast_n y a_n)\Delta_{23} (h)
 & = \Delta_{13} (h) (S(a_n)^\ast \otimes S(a_n) \otimes
y)\Delta_{23}(h)\\
& \rightarrow \Delta_{13} (h) (1 \otimes 1 \otimes y) \Delta_{23} (h).
\endsplit
$$
So we must have $y = (\iota \otimes \omega^\prime)\Delta(h)$.
\einspr

Having this result we can prove the main results of this section.

\inspr{5.8} Proposition \rm
  $\varphi$ is a faithful lower semi-continuous weight.

\snl\bf Proof : \rm
  It is clear from the definition 5.2 that $\varphi (\lambda x) = \lambda
\varphi (x)$ when $x \in A^\ast$ and $\lambda \geq 0$.  If $x,y \in A^+$ and
if $x = (\iota \otimes \omega ) \Delta (h)$ and $y = (\iota \otimes
\omega^\prime) \Delta (h)$ for some $\omega, \omega^\prime \in A^\ast_+$, then
$x+y = (\iota \otimes (\omega + \omega^\prime))\Delta (h)$ and we have $\|
(\omega + \omega^\prime) | J \| = \| \omega | J \| + \| \omega^\prime | J \|$ so
that $\varphi (x + y) = \varphi (x) + \varphi (y)$.  If $x,y \in A^\ast$ and
$\varphi(x) = \infty$ or $\varphi (y) = \infty$, then also $\varphi (x + y) =
\infty$ because when $x + y \in J^+_0$, also $x,y \in J^+_0$ by proposition 5.7.
\snl
We now prove that $\varphi$ is lower semi continuous.  Fix any $\lambda > 0$ and
consider a sequence $(x_n)$ in $A^+$ such that $x_n \rightarrow x$ and $\varphi
(x_n) \leq \lambda$.  Choose elements $\omega_n \in A^\ast_+$ such that $x_n =
(\iota \otimes \omega_n)\Delta(h)$.  Because $\varphi (x_n) \leq \lambda$ we get
$\| \omega_n | J \| \leq \lambda$.  By compactness we can assume that $\omega_n
\rightarrow \omega$ for some $\omega \in A^\ast_+$ in the
$\omega^\ast$-topology.  Then, as in 5.6, we get that $x = (\iota \otimes
\omega)\Delta(h)$. So also $\varphi (x) \leq \lambda$ because $\| \omega | J \|
\leq \lambda$.
\snl
If $x$ is positive and $\varphi(x)=0$, then $x=(\jo\ot\omega)(\Delta(h)$
for some $\omega\in A_+^\ast$. In this case, we must have $\omega| J=0$
and $x=0$. Hence, $\varphi$ is faithful.
\einspr

The weight need not be semi-finite because it is only finite on $J^+_0$ and this
need not be dense in $A^+$.  It is of course semi-finite in the discrete case.
\snl
Now we come to the invariance property.

\inspr{5.9} Proposition \rm
  If $a \in A^\ast$ and $\varphi(a) < \infty$ and if $\omega
\in A^\ast_+$, then
$$\varphi ((\iota \otimes \omega)\Delta(a)) = \| \omega | J \| \cdot
\varphi(a).$$

\snl\bf Proof : \rm
  Let $a = (\iota \otimes \omega^\prime)\Delta(h)$.  Then
$$\split
(\iota \otimes \omega) \Delta (a) & = (\iota \otimes \omega \otimes
\omega^\prime) \Delta^{(2)} (h) \\
& = (\iota \otimes \omega \omega^\prime) \Delta (h).
\endsplit
$$
So $\varphi((\jo\ot\omega)\Delta(a))=\|\omega\omega'|J\|$. If $(e_n)$ is
an approximate identity for $J$, then
   $$\|\omega\omega'|J\|=\lim_n(\omega\omega')(e_n)
   =\lim_n(\omega\ot\omega')\Delta(e_n).$$
Because $\Delta(J)(J\ot J)$ is dense in $J\overline\ot J$, we have that
still $(\Delta(e_n))_n$ is an approximate identity in $J\overline\ot J$.
Hence
  $$\lim_n(\omega\ot\omega')\Delta(e_n)
    =\|(\omega\ot\omega')|J\overline\ot J)\|
    =\|\omega|J\|\|\omega'|J\|.$$
Therefore
$$\split
\varphi ((\iota \otimes \omega) \Delta (a))
& = \| \omega \omega^\prime | J \|
\\ & = \| \omega | J \| \| \omega^\prime | J \| \\
& = \| \omega | J \| \cdot \varphi (a).
\endsplit
$$
\einspr

\nl

\bf 6. The regular representation \rm
\nl
Given the Haar weight $\varphi$, we can define the associated representation.
Let us recall this and see what we get here.
\snl
Denote, as usual, ${\frak N} = \{ a \in A \mid \varphi (a^\ast a) < \infty \}$.  So
$a \in {\frak N}$ iff $a^\ast a = (\iota \otimes \omega) \Delta (h)$ with $\omega
\in A^\ast_+$.  The set ${\frak M} = {\frak N}^\ast {\frak N}$ is in this case
precisely $J_0 = \{ (\iota \otimes \omega)\Delta (h) \mid \omega \in A^\ast \}$.
We have that ${\frak N}$ is a left ideal and that $\langle a,b \rangle
 = \varphi (b^\ast a)$
defines a scalar product on ${\frak N}$.  In this case, $\varphi$ is faithful and
so we do not have to divide by an ideal.  Denote by ${\Cal H}$ the
Hilbert space completion of $\frak N$ with respect to this scalar
product
and by $\Lambda (a)$ the image of an element $a$ in ${\frak N}$ in the Hilbert
space ${\Cal H}$.  The representation $\pi$ of $A$ on ${\Cal H}$ is given by
$\pi(a) \Lambda (b) = \Lambda (ab)$ when $a \in A$ and $b \in {\frak N}$.
Because the Haar weight is lower semi-continuous, this representation is
non-degenerate.  (See e.g. [16].)
\snl
We want to prove some specific results on this representation.

\inspr{6.1} Proposition \rm
  The space ${\Cal D}(S) \cap J_0$ is dense in ${\frak N}$ with
respect to the Hilbert space norm.

\snl\bf Proof : \rm
  We know already that
$$({\Cal D}(S) \cap J_0)^\ast ({\Cal D}(S) \cap J_0) \subseteq J_0$$
and so ${\Cal D(S)} \cap J_0 \subseteq {\frak N}$.  The more difficult
 part is to prove
that this is dense.
\snl
Take $x \in {\frak N}$ and $a \in {\Cal D}(S) \cap J_0$.
Then
 $xa \in {\frak N}$ because $\frak N$ is a left ideal.  If
$x^\ast x = (\iota \otimes \omega)\Delta(h)$, then
  $$\| \Lambda (x) - \Lambda (xa) \|^2
    = \varphi ((1 - a)^\ast x^\ast x(1-a))
    = \| \omega_a | J \|$$
where $\omega_a(z) = \omega ((1 - S(a)^\ast) z(1 - S(a)))$.
\snl
Now $\omega|J$ has the form $\omega_1 (\,\cdot \,p)$, where $p
\in J$
(see e.g. [25]).  We
know that $p$ can be approximated by elements $ep$ with $e \in J$ and $\| e \|
\leq 1$ (using an approximate identity).  Also $e$ can be approximated by
elements $S(a)$ with $a \in {\Cal D} (S) \cap J_0$ (see proposition 5.7). We can
assume $\| S(a) \| \leq 1$.  Hence we find $a \in {\Cal D} (S) \cap J_0$ such
that $\| (1-S(a))p\|$ is small and $\| S(a) \| \leq 1$.  Then $\|
\omega_a|J \|$
above is small.  This shows that $\frak N({\Cal D}(S) \cap J_0)$ is dense in
${\frak N}$.  In particular $J({\Cal D}(S) \cap J_0)$ is dense.
\snl
Finally, because ${\Cal D}(S) \cap J_0$ is norm dense in $J$, we have $({\Cal
D}(S) \cap J_0)({\Cal D}(S) \cap J_0)$ and hence ${\Cal D}(S) \cap J_0$
is dense in
${\frak N}$.
\einspr

In fact, in the argument above, we see that also the smaller set, namely the
elements of the form $(\iota \otimes \omega)\Delta(h)$ with $\omega(a) =
\psi(\Delta(h)(a \otimes 1))$ will still work.  Then $\omega \in {\Cal
D}(S_0^{-1})$.
\snl
Also remark that, not only $A$ acts non-degenerately, but this is also
true for $J$. This also follows from the arguments in the above proof.
\snl
In that case we have the following expression for the scalar product.

\inspr{6.2} Proposition \rm
  If $\omega_1,\omega_2 \in {\Cal D}(S_0^{-1})$ and $a = (\iota
\otimes \omega_1)\Delta(h)$ and $b = (\iota \otimes \omega_2)\Delta(h)$, then
$\varphi (b^\ast a) = (\omega_2^\ast \omega_1)(h)$.

\snl\bf Proof : \rm
  We have $a\in\Cal D(S)$ by proposition 4.8 and
$$
b^\ast a  = ((\iota \otimes \overline \omega_2)\Delta(h)) (a)
 = (\iota \otimes \overline \omega_2)(\Delta(h)(1 \otimes S(a)).
$$
So
$$\split
\varphi (b^\ast a)  = \overline \omega_2 (S(a))
& = \overline \omega_2 ((S_0^{-1} \omega_1 \otimes \iota)\Delta (h)) \\
& = (S_0^{-1} \omega_1 \otimes \overline \omega_2)\Delta(h) \\
& = (S_0^{-1} \omega_1 \otimes S_0^{-1} \omega^\ast_2)\Delta(h) \\
& = (\omega^\ast_2 \otimes \omega_1) \Delta (h) \\
& = (\omega^\ast_2 \omega_1)(h).
\endsplit
$$
\einspr

We will use $\Gamma(\omega)$ for $\Lambda((\iota \otimes \omega)\Delta(h))$.
Then we get a simple expression for the representation for the
representation of $A$.

\inspr{6.3} Proposition \rm
  If $a \in A$ then
$$\langle \pi(a) \Gamma (\omega_1), \Gamma(\omega_2) \rangle = (\omega_2^\ast \otimes
\omega_1)((a \otimes 1)\Delta(h)).$$

\snl\bf Proof : \rm
First, let $\omega_3 \in {\Cal D} (S_0^{-1})$
and
$a = (\iota \otimes \omega_3)\Delta(h)$ and denote $b = (\iota \otimes
\omega_1)\Delta(h)$.  Then $ab = (\iota \otimes \omega_3)(\Delta(h)(1 \otimes
S(b))$.  So
$$\langle \pi (a) \Gamma(\omega_1), \Gamma(\omega_2) \rangle = (\omega^\ast_2 \otimes
\omega_3) (\Delta (h)(1 \otimes S(b)).$$
Now, for all $x\in A$ we have
  $$\split
   \omega_3 (xS(b)) & = \omega_3 (x(S^{-1} \omega_1 \otimes \iota)\Delta(h)) \\
   & = (S^{-1} \omega_1 \otimes \omega_3)((1 \otimes x)\Delta(h)).
  \endsplit
  $$
If $x \in {\Cal D}(S)$ we have
 $$\split
 \omega_3(xS(b)) & = (S^{-1} \omega_1 \otimes \omega_3)((S(x) \otimes
1)\Delta(h)) \\
& = (S^{-1}  \omega_1)(S(x)a).
\endsplit
$$
Hence
  $$\langle\pi (a) \Gamma (\omega_1) , \Gamma (\omega_2)\rangle
    = (\omega_2^\ast \otimes
      S^{-1} \omega_1)((\iota \otimes S)\Delta(h) (1 \otimes a)).$$
By continuity, this holds for all $a\in J$.
Then, if now $a\in \Cal D(S^{-1})$ we get for the left hand side of the
above formula
  $$(\omega^\ast_2 \otimes \omega_1)(\iota \otimes S^{-1} (a)) \Delta(h))
    = (\omega^\ast_2 \otimes \omega_1)((a \otimes 1)\Delta (h)).
  $$
Again by continuity, this holds for all $a\in J$.
Because the representation on $J$ is already non-degenerate, this
formula also hods for $a\in A$.
\einspr

Using this formula, we can verify  that we have
a $^\ast$-representation. If $a\in\Cal D(S^{-1})$ we have
$$\split
\langle \Gamma (\omega_1), \pi(a^\ast) \Gamma(\omega_2) \rangle
& = \langle \pi(a^\ast)\Gamma(\omega_2),\Gamma(\omega_1) \rangle^- \\
& = (\omega^\ast_2 \otimes \omega_1)((a^\ast \otimes 1) \Delta (h))^- \\
& = (\omega^\ast_2 \otimes \omega^{\ast\ast}_1)(a^\ast \otimes 1)\Delta(h))^- \\
& = (\omega_2 \otimes \omega^\ast_1)((S \otimes S)(a^\ast \otimes
1)\Delta(h))^\ast) \\
& = (\omega_2 \otimes \omega^\ast_1) ((S^{-1}(a) \otimes 1) \sigma \Delta (h))
\\
& = (\omega^\ast_1 \otimes \omega_2)((a \otimes 1) \Delta(h)).
\endsplit
$$
We can also check positivity. Again, if $a\in\Cal D(S^{-1})$ we get
$$\split
(\omega^\ast \otimes \omega)((a^\ast a \otimes 1)\Delta (h)) & =
(\omega^\ast
\otimes \omega) (a^\ast \otimes S^{-1} (a)) \Delta (h)) \\
& =\omega^\ast (a^\ast \,\cdot\,) \cdot \omega (S^{-1} (a) \,\cdot\,)(h) \\
& =\omega (S(a^\ast)^\ast S(\,\cdot\,)^\ast)^- \omega (S^{-1}(a) \,\cdot\,)(h) \\
& =\omega (S^{-1} (a) S(\,\cdot\,)^\ast)^- \omega (S^{-1} (a) \,\cdot\,)(h).
\endsplit
$$
This is positive because it is of the form $\langle\Gamma(\omega'),
\Gamma(\omega')\rangle$ with $\omega'(z)=\omega(S^{-1}(a)z)$.
\snl
Now we want to
 define the associated representation of $A^\ast$
 by $$\pi(\omega)\Lambda(a) = \Lambda ((\iota \otimes \omega)\Delta(a)).$$
For this we need the following lemma.

\inspr{6.4} Lemma \rm
  If $a \in {\frak N}$ and $\omega \in A^\ast$, then $(\iota \otimes
\omega) \Delta (a) \in {\frak N}$.

\snl\bf Proof : \rm
  Let $\omega \in A^\ast$.  Consider $A$ in its universal representation.
Consider the polar decomposition of $\omega$ on $A^{\ast \ast}$.  So $\omega = |
\omega | (u \,\cdot\,)$ for some partial isometry.  Let $\xi$ be a vector such
that $| \omega | = \langle \,\cdot\, \xi , \xi \rangle$.
Then $\omega = \langle \,\cdot\, \xi, u^\ast \xi \rangle$.
\snl
Let $a \in {\frak N}$ and put $b = (\iota \otimes \omega)\Delta(a)$.  Let $\eta$
be any vector and $(e_i)_{i \in I}$ a basis in the Hilbert space on which $A$
acts.  Then
$$\split
\langle b^\ast b \eta, \eta \rangle
 = \sum_i | \langle b \eta, e_i \rangle |^2
& = \sum_i | \langle \Delta(a) \eta \otimes \xi, e_i \otimes u^\ast \xi \rangle |^2 \\
& \leq \| u^\ast \xi \|^2 \| \Delta (a) \eta \otimes \xi \|^2 \\
& \leq \| \xi \|^2 \langle \Delta (a^\ast a) \eta \otimes \xi, \eta \otimes \xi \rangle \\
& = \| \xi \|^2 \langle (\iota \otimes | \omega |) \Delta (a^\ast a) \eta, \eta \rangle.
\endsplit
$$
So
$$b^\ast b \leq \| \omega \| (\iota \otimes | \omega | ) \Delta (a^\ast a).$$
We know that $a^\ast a \in J_0$ and so also $(\iota \otimes | \omega
|)\Delta(a^\ast a) \in J_0$.
\snl
Then $b^\ast b \in J_0$ and $b \in {\frak N}$.
\einspr

If we apply $\varphi$ to the above equation, we get
  $$\varphi (b^\ast b)  \leq \| \omega \| \varphi ((\iota \otimes | \omega |) \Delta
    (a^\ast a))  = \| \omega \| \| \omega \| \varphi (a^\ast a).$$
So we can define a bounded operator on ${\Cal H}$.

\inspr{6.5} Definition \rm
  If $\omega \in A^\ast$ we define $\pi(\omega)$ on ${\Cal H}$
by $\pi (\omega) \Lambda (a) = \Lambda (b)$ where $a \in {\frak N}$ and $b =
(\iota \otimes \omega) \Delta (a)$.
\einspr

We see that $\| \pi (\omega) \| \leq \| \omega \|$.  We also have $\pi
(\omega_1 \omega_2) = \pi (\omega_1) \pi (\omega_2)$ when $\omega_1, \omega_2
\in A^\ast$ and $\omega_1 \omega_2$ is defined as before by $(\omega_1
\omega_2)(a) = (\omega_1 \otimes \omega_2) \Delta (a)$.  We also have that $\pi
(\omega^\ast) = \pi (\omega)^\ast$ when $\omega \in {\Cal D} (S_0^{-1})$.  This
will follow from the following lemma.

\inspr{6.6} Lemma \rm
  If $\omega,\omega_1 \in {\Cal D} (S_0^{-1})$, then
$\pi (\omega) \Gamma (\omega_1) = \Gamma (\omega \omega_1).$

\snl\bf Proof : \rm
  Recall that $\Gamma (\omega_1) = \Lambda (a)$ when $a = (\iota \otimes
\omega_1) \Delta (h)$.  So
$\pi (\omega) \Gamma (\omega_1) = \Lambda (b)$
where $b = (\iota \otimes \omega) \Delta (\iota \otimes \omega_1)\Delta (h)) =
(\iota \otimes \omega \omega_1) \Delta (h)$.  Hence $\pi (\omega) \Gamma
(\omega_1) = \Gamma (\omega \omega_1)$.
\einspr

We have seen that $\langle \Gamma (\omega_1), \Gamma (\omega_2) \rangle
= (\omega^\ast_2
\omega_1) (h)$.  So $\langle \pi (\omega) \Gamma (\omega_1), \Gamma (\omega_2) \rangle =
(\omega^\ast_2 (\omega \omega_1)) (h) = ((\omega^\ast \omega_2)^\ast
\omega_1)(h) = \langle \Gamma (\omega_1), \pi (\omega^\ast) \Gamma (\omega_2) \rangle$.  So $\pi$
is a $^\ast$-representation.
\snl
This specific combination of the representation of $A$ and of $A^\ast$
will give us the fundamental unitary in the next section.

\nl\nl

\bf 7. The fundamental unitary \rm
\nl
We have a representation of $A$ and $A^\ast$ on ${\Cal H}$.  This gives us a
representation of $A^\ast \otimes A$ on $\Cal H\overline\ot \Cal H$.  If $\omega \in A^\ast$ and $a \in A$ we
can write
$$\split
 \langle \pi(\omega) \otimes \pi(a) \xi_1 \otimes \eta_1, \xi_2 \otimes \eta_2 \rangle
& = \langle \pi (\omega) \xi_1, \xi_2 \rangle \langle \pi (a) \eta_1 , \eta_2 \rangle \\
& = \langle \pi (\gamma) \xi_1, \xi_2 \rangle
\endsplit
$$
where $\gamma(b) = \langle \pi (a) \eta_1, \eta_2 \rangle \omega (b) = \langle \pi (w(b)a) \eta_1,
\eta_2 \rangle$ for all $b\in A$. If we consider $\omega \otimes a$ as the linear map $f
: b \rightarrow \omega (b)a$ from $A$ to $A$ we find
$$\gamma = \omega_{\eta_1,\eta_2} \circ \pi \circ f.$$
With $f = \iota$ we will obtain the fundamental operator $W$.  So, we should
have
$$\langle W \xi_1 \otimes \eta_1, \xi_2 \otimes \eta_2 \rangle = \langle \pi
(\omega_{\eta_1,\eta_2}) \xi_1, \xi_2 \rangle.$$
We first prove that we can define a bounded operator by this formula.

\inspr{7.1} Proposition \rm
  There is a bounded operator $W$ on ${\Cal H} \otimes {\Cal
H}$ defined by
$$\langle W \xi_1 \otimes \eta_1, \xi_2 \otimes \eta_2 \rangle = \langle \pi
(\omega_{\eta_1,\eta_2}) \xi_1, \xi_2 \rangle.$$

\snl\bf Proof : \rm
  Take a vector $\sum_k \xi_k \otimes \eta_k$ in the algebraic tensor
product of ${\Cal H}$ with itself.  Let $(e_i)_i$ be an orthonormal basis in
${\Cal H}$.
\snl
Consider the numbers $(p_{ij})_{i,j}$ defined by
$$p_{ij} = \sum_k \langle \pi (\omega_{\eta_k,e_i}) \xi_k, e_j \rangle.$$
Then
$$\split
\sum_j | p_{ij} |^2 & = \sum_{k,\ell,j} \langle \pi (\omega_{\eta_\ell,
e_i}) \xi_\ell,
e_j \rangle^- \langle \pi (\omega_{\eta_k,e_i}) \xi_k, e_j \rangle \\
& = \sum_{k,\ell} \langle \pi (\omega_{\eta_k,e_i}) \xi_k, \pi
(\omega_{\eta_\ell,e_i}) \xi_\ell \rangle.
\endsplit
$$
Now assume that $\xi_k = \Lambda (a_k)$ with $a_k \in {\frak N}$.  Then
$$\split
 \langle \pi (\omega_{n_k,e_i}) \xi_k, \pi (\omega_{\eta_\ell,e_i}) \xi_\ell \rangle
& = \langle \Lambda ((\iota \otimes \omega_{\eta_k,e_i}) \Delta (a_k)), \Lambda ((\iota
\otimes \omega_{\eta_\ell,e_i}) \Delta (a_\ell)) \rangle \\
& =  \varphi (((\iota \otimes \omega_{\eta_\ell,e_i}) \Delta (a_\ell))^\ast
(\iota \otimes \omega_{\eta_k,e_i}) \Delta (a_k)).
\endsplit
$$
If $p_i$ is the projection on the one-dimensional subspace $\Bbb C e_i$, we
can rewrite this as
$$\varphi ((\iota \otimes \omega_{\eta_k,\eta_\ell})(\Delta(a^\ast_\ell)(1
\otimes p_i) \Delta (a_k))).$$
Any finite sum over $i$ will remain smaller than
$$ \varphi (\sum_{k,\ell} (\iota \otimes \omega_{\eta_k, \eta_\ell}) \Delta
(a^\ast_\ell a_k))
 = \sum _{k,\ell}\varphi (a^\ast_\ell a_k) \omega_{\eta_k, \eta_\ell} (1)
 = \| \sum_k
\xi_k \otimes \eta_k \|^2.
$$
\einspr

This shows also that $\| W \| \leq 1$.  In the following proposition, we show that
$W^\ast W = 1$.

\inspr{7.2} Proposition \rm
 $W^\ast W = 1$.

\snl\bf Proof : \rm
  When we look at the proof of the previous proposition, we must show that
$$\sum_i \sum_{k,\ell} (\iota \otimes \omega_{\eta_k, \eta_\ell})(\Delta
(a^\ast_\ell)(1 \otimes p_i)\Delta(a_k))$$
converges in norm to
$$\sum_{k,\ell} (\iota \otimes \omega_{\eta_k,\eta_\ell}) \Delta (a^\ast_\ell
a_k).$$
Then we can use the lower semi-continuity of $\varphi$.
\snl
We know that $\sum p_i = 1$ in the weak operator topology.  Therefore, it is
sufficient to show that the map
$$x \rightarrow (\iota \otimes \omega_{\xi, \eta})(\Delta(b) (1 \otimes
x)\Delta(c))$$
is continuous from $\Cal B({\Cal H})$ with the weak operator topology to $A$ with the
norm topology.  We can work on the unit ball of $\Cal B({\Cal H})$.
\snl
If $\xi = p \xi_1$ and $\eta = q \eta_1$ with $p,q \in A$, we have
$$
(\iota \otimes \omega_{\xi,\eta})  (\Delta (b)(1 \otimes x)\Delta(c))
= (\iota \otimes \omega_{\xi_1,\eta_1})((1 \otimes q^\ast)\Delta(b)
 (1 \otimes x) \Delta (c) (1 \otimes p)).
$$
Now $\Delta (c)(1 \otimes p)$ and $(1 \otimes q^\ast)\Delta(b)$ can be
approximated in norm by elements in the algebraic tensor product of $A$ with
itself.  But clearly, if $p_1,q_1,p_2,q_2 \in A$, then
$$x \rightarrow (\iota \otimes \omega_{\xi_1,\eta_1})((p_1 \otimes q_1)(1
\otimes x)(p_2 \otimes q_2))$$
has the correct continuity property.  And since we work on the unit ball of
${\Cal B}({\Cal H})$, this will remain true in the limit.
\einspr

We can not use the same argument to show that also $WW^\ast = 1$ so that $W$ is
a unitary.  This is of a different nature.  The properties of the antipode are
needed here.  We will first prove the other properties of $W$.  Then the
unitarity will follow easily.
\nl
In the following proposition, we will assume that $A$ acts on ${\Cal H}$.  Then
$M(A \overline\otimes A)$ acts on ${\Cal H} \overline{\otimes} {\Cal H}$.
So we can avoid the use of too many $\pi$'s.

\inspr{7.3} Proposition \rm
  For any $a \in A$ we get
$W(a \otimes 1)  = \Delta(a)W.$

\snl\bf Proof : \rm  Take $\omega_1, \omega_2 \in A^\ast$
and $a \in {\Cal D} (S^{-1})$, $\omega_2 \in {\Cal D} (S^{-1})$.  Then
$$\split
 \langle W(a \Gamma (\omega_1) \otimes \xi) \rangle , \Gamma (\omega_2)
 \otimes \eta \rangle
& = \langle \pi (\omega_{\xi,\eta}) a \Gamma (\omega_1), \Gamma (\omega_2)
\rangle \\
& = \langle \pi (\omega_{\xi,\eta}) \Gamma (\omega_1 (S^{-1} (a) \cdot )),
\Gamma
(\omega_2) \rangle \\
& = (\omega_2^\ast \otimes \omega_{\xi, \eta} \otimes \omega_1 (S^{-1} (a) \cdot))
\Delta^{(2)} (h) \\
& = ((\omega_2^\ast \omega_{\xi, \eta}) \otimes \omega_1)(1 \otimes
S^{-1}(a))\Delta(h))\\
& = ((\omega_2^\ast \omega_{\xi, \eta}) \otimes \omega_1)((a \otimes 1) \Delta (h))
\\
& = (\omega_2^\ast \otimes \omega_{\xi, \eta} \otimes \omega_1)((\Delta(a)
\otimes 1) \Delta^2 (h)).
\endsplit
$$
On the other hand, if $p^\ast \in {\Cal D} (S^{-1})$ and $q \in A$, then
$$\split
\langle (p \otimes q)W(\Gamma(\omega_1) \otimes \xi), \Gamma(\omega_2) \otimes \eta \rangle
& = \langle W(\Gamma(\omega_1) \otimes \xi), p^\ast \Gamma (\omega_2) \otimes q^\ast
\eta \rangle \\
& = \langle \pi (\omega_{\xi,q^\ast \eta}) \Gamma (\omega_1), \Gamma (\omega_2 (S^{-1}
(p^\ast) \,\cdot\,)) \rangle \\
& = ((\omega_2 (S^{-1} (p^\ast) \,\cdot\,))^\ast \otimes \omega_{\xi,q^\ast \eta} \otimes
\omega_1) (\Delta^{(2)} (h)).
\endsplit
$$
Now
$$\split
(\omega_2 (S^{-1} (p^\ast) \cdot )^\ast) (a) & = \omega_2 (S(p)^\ast S(a)^\ast)^- \\
& = \omega_2 (S(pa)^\ast)^- = \omega^\ast_2 (pa).
\endsplit
$$
So
$$\split
\langle (p \otimes q)W(\Gamma(\omega_1) \otimes \xi), \Gamma (\omega_2) \otimes \eta \rangle
& = (\omega^\ast_2 (p \,\cdot\,) \otimes \omega_{\xi, q^\ast \eta} \otimes
\omega_1)(\Delta^{(2)} (h)) \\
& = (\omega^\ast_2 \otimes \omega_{\xi,\eta} \otimes \omega_1)((p \otimes q
\otimes 1) \Delta^{(2)} (h)).
\endsplit
$$
We can take $\eta = c\eta_1$ and approximate $(1 \otimes c^\ast) \Delta (a)$
by linear combinations of $p \otimes q$. Then, by continuity we have
$$\split
\langle \Delta (a) W (\Gamma (\omega_1) \otimes \xi), \Gamma (\omega_2) \otimes \eta \rangle &
= (\omega^\ast_2 \otimes \omega_{\xi, \eta} \otimes \omega_1)((\Delta(a) \otimes
1) \Delta^{(2)} (h)) \\
& = \langle W (a \Gamma (\omega_1) \otimes \xi), \Gamma (\omega_2) \otimes \eta \rangle
\endsplit
$$
\einspr

We now want to prove $(\iota \otimes \Delta)W = W_{12} W_{13}$.  We can give a
meaning to this formula by applying $\omega \otimes \iota \otimes \iota$ for
some $\omega \in {\Cal B} (\Cal H)_\ast$.  Indeed, $(\omega \otimes \iota)(W)
\in A$.  To show this, consider $\langle W \Gamma (\omega_1) \otimes \xi, \Gamma
(\omega_2) \otimes \eta \rangle$.  This is equal to $(\omega^\ast_2 \otimes
\omega_{\xi,\eta} \otimes \omega_1)(\Delta^2 (h)) = \langle a\xi, \eta \rangle$ where $a
= (\omega^\ast_2 \otimes \iota \otimes \omega_1) (\Delta^2 (h))$.

\inspr{7.4} Proposition \rm
  $(\iota \otimes \Delta) W = W_{12} W_{13}$.

\snl\bf Proof : \rm
  By the remark above, we have to show $\Delta ((\omega \otimes \iota)W) =
(\omega \otimes \iota \otimes \iota)(W_{12} W_{13})$.  We will do this for
$\omega = \omega_{\xi,\eta}$.  We have for $\xi, \eta, \xi_1, \eta_1,  \xi_2,
\eta_2 \in {\Cal H}$,
$$\split
\langle \Delta ((\omega_{\xi,\eta} \otimes \iota)W) \xi_1 \otimes \eta_1, \xi_2
\otimes \eta_2 \rangle & = \langle \pi(\omega_{\xi_1,
 \xi_2} \omega_{\eta_1, \eta_2}) \xi,
\eta \rangle \\
& = \langle \pi (\omega_{\xi_1, \xi_2})
 \pi (\omega_{\eta_1, \eta_2}) \xi, \eta \rangle \\
& = \langle W \pi (\omega_{\eta_1 , \eta_2}) \xi \otimes \xi_1 , \eta \otimes \xi_2 \rangle
\\
& = \langle \pi (\omega_{\eta_1 ,\eta_2}) \xi \otimes \xi_1, W^\ast (\eta \otimes
\xi_2) \rangle \\
& = \langle W_{13} \xi \otimes \xi_1 \otimes \eta_1 , W^\ast_{12} (\eta \otimes \xi_2
\otimes \eta_2) \rangle \\
& = \langle W_{12} W_{13} \xi \otimes \xi_1 \otimes \eta_1 , \eta
\otimes \xi_2
\otimes \eta_2 \rangle.
\endsplit
$$

\inspr{7.5} Proposition \rm
  $W$ satisfies the Pentagon equation
$$W_{23} W_{12} = W_{12} W_{13} W_{23}.$$

\snl\bf Proof : \rm  Applying $\omega \otimes \iota \otimes \iota$, we see that we need to
prove, when $a = (\omega \otimes \iota)(W)$ that
$$W(a \otimes 1) = (\omega \otimes \iota \otimes \iota)(W_{12} W_{13})W.$$
But we have shown in 7.3 that
$$W(a \otimes 1) = \Delta (a) W$$
and in 7.4 that
$$\Delta (a) = (\omega \otimes \iota \otimes \iota)(W_{12} W_{13}).$$
Combining these results, we find the Pentagon equation.
\einspr

Now, it is not so difficult anymore to show that $W$ is a unitary.

\inspr{7.6} Proposition \rm
  $W$ is unitary.

\snl\bf Proof : \rm  Apply $\iota \otimes \iota \otimes \omega$ to the Pentagon equation.
Then
$$(1 \otimes b)W = W((\iota \otimes \iota \otimes \omega)(W_{13}
W_{23}))$$
where $b = (\iota \otimes \omega) W = \pi (\omega)$.
We see that $1 \otimes \pi(\omega)W$ maps the range of $W$ into itself.
Now
$$\langle W (\xi \otimes \Lambda (h)), \eta \otimes \gamma  \rangle
=\langle \pi(\omega_{\Lambda (h), \gamma}) \xi, \eta \rangle.$$
If $\xi = \Lambda (a)$ then
$$\pi (\omega_{\Lambda (h), \gamma}) \xi = \Lambda ((\iota \otimes
\omega_{\Lambda (h), \gamma}) \Delta (a)).$$
But $\Lambda (h) = h \Lambda (h)$, so we get
$$(\Lambda (\iota \otimes \omega_{\Lambda (h)), \gamma})(\Delta (a) (1 \otimes h))
= \omega_{\Lambda (h), \gamma} (1) \Lambda (a).$$
Hence $\pi(\omega_{\Lambda (h), \gamma}) = \omega_{\lambda (h), \gamma} (1) 1$.
And $\langle W \xi \otimes \Lambda (h), \eta \otimes \gamma \rangle = \langle \xi, \eta \rangle \langle
\Lambda (h), \gamma \rangle$.  So $W \xi \otimes \Lambda (h) = \xi \otimes \Lambda
(h)$.  Then, also $(1 \otimes \pi(\omega) (\xi \otimes \Lambda(h))$ is in the
range of $W$, and this is precisely $\xi \otimes \Gamma (\omega)$.  These
vectors span ${\Cal H} \overline\otimes {\Cal H}$.
\einspr

Having this fundamental operator satisfying the Pentagon equation, one
obtains the (reduced) dual by standard methods.
\nl\nl

\noindent
\bf  8. References \rm
\bigskip\parindent=20pt
\item{[1]} E.\ Abe: \it Hopf Algebras. \rm Cambridge University Press (1977).

\item{[2]} S.\ Baaj: {\it R\'epresentation r\'eguli\`ere du groupe quantique
         $E_{\mu}(2)$ de Woronowicz}. C.R.\ Acad.\ Sci., Paris, {\bf 314} (1992), 1021--1026.
 
\item{[3]} S.\ Baaj: {\it  R\'epresentation r\'eguli\`ere du groupe quantique
         des d\'eplacements de Woronowicz}. Recent advances in operator algebras (Orl\'eans, 1992), Ast\'erisque {\bf 232} (1995), 11--48.

\item{[4]} S.\ Baaj \& G.\ Skandalis: \it Unitaires Multiplicatifs et dualit\'e
          pour les Produits Crois\'es de C$^*$-Alg\`ebres. \rm Ann.\ Scient.\ Ec.\ Norm.\ Sup., {\bf 26} (1993), 425--488.

\item{[5]} V.G.\ Drinfel'd: \it Quantum Groups. \rm Proceedings ICM Berkeley
          (1986), 798--820.

\item{[6]} E.G.\ Effros \& Z.-J.\ Ruan: {\it Discrete Quantum Groups I.
        The Haar Measure}. Int.\ J.\ Math.\ {\bf 5} (1994), 681--723.

\item{[7]} M.\ Enock \& J.-M.\ Schwartz: \it Kac Algebras and Duality of
         Locally Compact Groups. \rm Springer, Berlin (1992).

\item{[8]} K.H.\ Hofmann: {\it Elements of compact semi-groups.} Charles
E. Merill Books Inc. Columbus, Ohio (1996).

\item{[9]} J.\ Hollevoet \& J.\ Quaegebeur: {\it Duality for Quantum
        Groups}. Preprint K.U.Leuven (1992). Unpublished.

\item{[10]} J.\ Kustermans: {\it Quasi-discrete quantum groups are almost discrete.} Preprint K.U.Leu\-ven (1995). Unpublished.

\item{[11]} J.\ Kustermans \& S.\ Vaes: {\it A simple definition for locally compact quantum groups.} C.R.\ Acad.\
Sci., Paris, S\'er.\ I Math.\ {\bf 328} (1999), 871--876.

\item{[12]} J.\ Kustermans \& S.\ Vaes: {\it Locally compact quantum groups}. Ann.\ Scient \'{E}c.\ Norm.\ Sup.\ {\bf 33} (2000), 837--934.

\item{[13]} J.\ Kustermans \& S.\ Vaes: {\it Locally compact quantum groups in the von Neumann algebraic setting.}
Math.\ Scand.\ {\bf 92} (2003), 68--92.

\item{[14]} M.B.\ Landstad \& A.\ Van Daele: {\it Compact and discrete subgroups of algebraic quantum groups.} Preprint University of Trondheim and K.U.Leuven (2004).  

\item{[15]} B.D.\ Malviya \& B.J.\ Tomiuk: {\it Multiplier operators on
B$^*$-algebras.} Proc.\ Amer.\ Math.\ Soc.\ {\bf 31} (1972) 505--510.

\item{[16]} G.K.\ Pedersen: \it C$^*$-algebras and their automorphism
groups. \rm Academic Press, London (1979).

\item{[17]} P.\ Podle\'s \& S.L.\ Woronowicz: \it Quantum deformation of Lorentz
          group. \rm Commun.\ Math.\ Phys.\ \bf 130 \rm (1990), 381--431.

\item{[18]} S.\ Sakai: \it C$^*$-algebras and W$^*$-algebras. \rm Springer
        Verlag (Berlin) 1971.

 \item{[19]} G.\ Skandalis  \it Operator Algebras and Duality. \rm Proceedings
           ICM Kyoto (1990), 997--1009.

 \item{[20]} M.E.\ Sweedler: \it Hopf Algebras. \rm Mathematical Lecture Note
          Series. Benjamin. N.Y., 1969.

\item{[21]} D.C.\ Taylor: {\it The strict topology for double centralizers
algebras}. Trans.\ Amer.\ Math.\ Soc.\ {\bf 150} (1970) 633--643.

\item{[22]} J.-M.\ Vallin: {\it C$^*$-alg\`ebres de Hopf et
       C$^*$-algebr\`es de
          Kac}.  Proc.\ Lond.\ Math.\ Soc.\ {\bf 50}  (1965)  131--174.

\item{[23]} A.\ Van Daele: {\it Dual pairs of Hopf $^*$-algebras.} Preprint K.U.Leuven (1990) (first version); Bull.\ Lond.\ Math.\ Soc.\ {\bf
          25} (1993) 209--230 (second version).

\item{[24]} A.\ Van Daele: \it Multiplier Hopf Algebras. \rm
          Trans.\ Amer.\ Math.\ Soc.\ {\bf 342} (1994), 917--932.

\item{[25]} A.\ Van Daele: {\it The Haar measure on a compact quantum
group.} Proc.\ Amer.\ Math.\ Soc.\ {\bf 123} (1995), 3125--3128.

\item{[26]} A.\ Van Daele: {\it Discrete quantum groups.} J.\ of Alg.\ {\bf 180} (1996), 431--444.

\item{[27]} A.\ Van Daele: {\it The Haar measure on finite quantum groups.} Proc.\ Amer.\ Math.\ Soc.\ {\bf 125} (1997), 3489--3500.

\item{[28]} A.\ Van Daele \& S.\ Wang: {\it The Larson-Sweedler theorem for multiplier Hopf algebras.} Preprint K.U.Leuven (2004).  
 
\item{[29]} S.L.\ Woronowicz: {\it Pseudospaces, pseudogroups and
Pontryagin duality.} Proceedings of the International Conference on
Mathematics and Physics, Lausanne (1979), Lecture Notes in Physics {\bf
116} () 407--412.

\item{[30]} S.L.\ Woronowicz: {\it Compact Matrix Pseudogroups.} Comm.\ Math.\
          Physics {\bf 111} (1987), 613--665.

\item{[31]} S.L.\ Woronowicz: {\it Compact quantum groups.}  Sym\'etries quantiques (Les Houches, 1995), 845--884, North Holland, Amsterdam (1998).

\end